\documentclass[12pt]{iopart}
\usepackage[latin1]{inputenc}
\usepackage[english]{babel}
\usepackage{iopams, setstack,enumerate}
\newcommand{\R}{\mathbb{R}}
\newcommand{\C}{\mathbb{C}}
\newcommand{\N}{\mathbb{N}}
\newcommand{\Z}{\mathbb{Z}}

\newcommand{\re}{\mathrm{Re}\,}
\newcommand{\im}{\mathrm{Im}\,}
\newcommand{\indic}{\mathrm{l\negthinspace l}}
\newcommand{\id}{\mathrm{Id}}
\newcommand{\ja}[2]{j_{#1}\left(\lambda  {#2}\right)}
\newcommand{\na}[2]{\eta_{#1}\left(\lambda  {#2}\right)}
\newcommand{\n}[2]{{\left\|{#1}\right\|}_{#2}}
\newcommand{\ing}{\lshad}
\newcommand{\ind}{\rshad}
\newcommand{\vect}[1]{\mathrm{Span}\left\{#1\right\}}
\newcommand{\sgn}{\mathrm{sgn}}
\newcommand{\G}{\mathcal{G}}
\newcommand{\classe}[1]{\mathcal{C}^{#1}}
\newcommand{\hunr}{H^1_\R(0,1)}
\newcommand{\hunc}{H^1_\C(0,1)}
\newcommand{\lc}{L_{\C}^2(0,1)}
\newcommand{\lr}{L_{\R}^2(0,1)}

\newcommand{\sr}{\mathcal{R}}
\renewcommand{\ss}{\mathcal{S}}

    \newtheorem{thm}{Theorem}[section]
    \newtheorem{coro}{Corollary}[section]
    \newtheorem{prop}{Proposition}[section]
    \newtheorem{lem}{Lemma}[section]
    \newtheorem{lema}{Lemma}[section]
    \newtheorem{defn}{Definition}[section]
    \newtheorem{nota}{Notations}
\newcommand{\qed}{$\square$}
\newcommand{\eqref}[1]{\eref{#1}}
\renewcommand{\H}{\mathcal{H}}
\renewcommand{\wr}{\mathcal{W}}
\renewcommand{\O}[1]{\mathcal{O}\negthinspace\left(#1\right)}

\begin{document}
%%%%%%%%%%%%%%%%%%%%%%%%%%%%%%%%%%%%%%%%%%%%%%%%%%%%%%%%%%%%%%%%%%%%%%%%%%%%%%%%%%%%%%%%

\title{Inverse spectral problem for singular AKNS operator on $[0,1]$.}
\author{Frédéric SERIER\footnote{Present address:
Institut für Mathematik, Universität Zürich, Winterthurerstrasse
190 CH-8057 Zürich}}
\address{LMJL - Laboratoire de Mathématiques Jean Leray\\ UMR CNRS 6629-UFR Sciences et
Techniques\\
2 rue de la Houssinière - BP 92208\\ F-44322 Nantes Cedex 3}
\ead{\mailto{frederic.serier@univ-nantes.fr}}
\begin{abstract}
We consider an inverse spectral problem for a class of singular
AKNS operators $H_a, a\in\N$ with an explicit singularity. We
construct for each $a\in\N$, a standard map
$\lambda^a\times\kappa^a$ with spectral data $\lambda^a$ and some
norming constant $\kappa^a$. For $a=0$, $\lambda^a\times\kappa^a$
was known to be a local coordinate system on $\lr\times\lr$. Using
adapted transformation operators, we extend this result to any
non-negative integer $a$, give a description of isospectral sets
and we obtain a Borg-Levinson type theorem.
\end{abstract}
\submitto{\IP} \maketitle

%   Structure principale
%\pagestyle{headings}
     \section{Introduction}

%   Introduction
The Schrödinger  operator $\H=-\Delta+q(\n{x}{})$ with a radial
potential $q$, acting on the unit ball of $\R^3$, through a
decomposition via spherical harmonics (see \cite{rs}, p.
$160-161$), is unitary equivalent to a collection of singular
differential operators $\H_a(q)$, $a\in \N$ acting on $\lr$, with
Dirichlet boundary conditions, defined by
\begin{equation*}
\fl\qquad\H_a(y)(x):=   \left(-\frac{d^2}{d
x^2}+\frac{a(a+1)}{x^2}+q(x)
    \right)y(x)=\lambda y(x),\quad x\in[0,1], \lambda\in\C.
\end{equation*}
With this splitting, it makes sense to study inverse spectral
problems not for $\H$ itself but for each $\H_a$.

The inverse spectral problem for these operator is the
construction for each $a\in\N$, of a regular coordinate system
$\lambda^a\times\kappa^a$ for potentials $q\in\lr$ where
$\lambda^a$ represent the spectrum of $\H_a$ and $\kappa^a$ are
convenient complementary data (regularity means stability of the
inverse spectral problem).

This question is not new and has been answered: Borg \cite{bo} and
Levinson \cite{lev} first, proved that $\lambda^0\times\kappa^0$
was one-to-one on $\lr$; then Pöschel and Trubowitz \cite{ist}
completed this result obtaining $\lambda^0\times\kappa^0$ as a
global real-analytic coordinate system on $\lr$. Guillot and
Ralston \cite{gr} extended their results to
$\lambda^1\times\kappa^1$, passing through the singularity inside
the equation. Next Zhornitskaya and Serov \cite{zs}, and Carlson
\cite{carlson}, proved that for all real $a\geq -1/2$,
$\lambda^a\times\kappa^a$ is one-to-one on $\lr$. Finally, the
author \cite{moi-SL} completed theses works proving that for all
$a\in\N$ the map $\lambda^a \times\kappa^a$ was a local (hence
global) diffeomorphism on $\lr$.

Then, it is natural and interesting to wonder if these kind of
results can be found for an other physical equation: the Dirac
equation. Hence, as the radial Schrödinger operator, the Dirac
operator with a radial electric potential acting on the unit ball
of $\R^3$ is decomposed (see for instance \cite{tha}) into a
collection of operators $H_a$ defined on $[0,1]$ by
\begin{equation}\label{akns-AKNS+V}
\fl \qquad H_a(V)Y(x):=\Bigg(\Bigg[\begin{array}{cc}0&-1\\
1&0\end{array}\Bigg]\frac{d}{d x}+\Bigg[\begin{array}{cc}0&-\frac{a}{x}\\
-\frac{a}{x}&0\end{array}\Bigg]+V(x)\Bigg)Y(x)=\lambda Y(x),
\end{equation}
where $Y=(Y_1,Y_2)$, $\lambda\in\C$ and
$$V(x)=\Bigg[\begin{array}{cc}q(x)+m&0\\0&q(x)-m\end{array}\Bigg], m\in\R;$$
with general boundary conditions
\begin{equation}\label{akns-AKNS-bord}
    Y_2(0)=0;\quad Y(1)\cdot u_\beta=0 \quad %\mathrm{where}
    \quad u_\beta=\Bigg[\begin{array}{c}
      \sin{\beta} \\
      \cos{\beta} \\
    \end{array}\Bigg],\quad \beta\in\R.
\end{equation}

Written this way, the Dirac operator seems to be unadapted in view
of inverse spectral problems. Indeed for $a=0$, as raised by
Levitan and Sargsjan in \cite{lesa-dirac}(Chap. 7) and pointed out
more generally by Clark and Gesztesy in \cite{cg} (section 6), the
existence of a gauge transformation on the potential $V$ leaving
the spectrum invariant leads to choose a normal form for the
problem, namely, the AKNS system, obtained from
\eqref{akns-AKNS+V} considering potentials $V$ of the following
shape:
\begin{equation}\label{akns-V_akns}
    V(x)=\Bigg[\begin{array}{cc}-q(x)&p(x)\\p(x)&q(x)\end{array}\Bigg],\quad(p,q)\in\lr\times\lr.
\end{equation}

Moreover, there are some clues showing that the inverse spectral
problem is kind of degenerated: for instance, the Ambarzumian type
theorem obtained by Kiss \cite{kiss} who proves that for all
$m\neq 0$ and $q\in\mathcal{C}([0,1];\R)$, if $H_0(V)$ has the
same eigenvalues as $H_0(0)$ then $q=0$. An other reason, to turn
to the AKNS operator, is it similarity with the Schrödinger
operator as figured out the papers of Grébert and Guillot
\cite{gg} and Amour and Guillot \cite{ag}. And finally, technical
difficulties arise when computing asymptotics for solutions of the
Dirac equation, see remark page \pageref{remarque-dirac}.
%For details and partial results, see \cite{moi-these}.

%%%%%%%%%%%%%%%%%%%%%%%%%%%%%%%%%%%%%%%%%%%%%%%%%%%%%%%%%%%%%%%%%%%%%%%%%%%
%%%%%%%%%%%%%%%%%%%%%%%%%%%%%%%%%%%%%%%%%%%%%%%%%%%%%%%%%%%%%%%%%%%%%%%%%%%

Our purpose is the stability of the inverse spectral problem for
$H_a$
(\eqref{akns-AKNS+V}-\eqref{akns-AKNS-bord}-\eqref{akns-V_akns}).
For this, we construct for each $a\in\N$, a spectral map
$\lambda^a\times\kappa^a$ for potentials $V$ with spectral data
$\lambda^a$ and some norming constant $\kappa^a$. The framework is
the work of Grébert and Guillot \cite{gg} for the regular operator
($a=0$). They constructed a local coordinate system
$\lambda^0\times\kappa^0$ on $\lr\times\lr$ and proved it is
global on $H_\R^j(0,1)\times H_\R^j(0,1)$ for $j=1,2$. With the
singularity, interesting problem arise and add supplementary
difficulties, especially when we study the invertibility of the
Fréchet derivative of $\lambda^a\times\kappa^a$. For this, we use
some transformation operators who, roughly speaking, reduce the
singularity.

Our result is that for all $a\in\N$ , $\lambda^a\times\kappa^a$ is
a local diffeomorphism on $\lr\times\lr$ and one-to-one on
$H^1_\R(0,1)\times H^1_\R(0,1)$. Moreover, we locally describe
sets of isospectral potentials as smooth submanifolds of
$\lr\times\lr$ with explicitly tangent and normal spaces.

%%%%%%%%%%%%%%%%%%%%%%%%%%%%%%%%%%%%%%%%%%%%%%%%%%%%%%%%%%%%%%%%%%%%%%%%%%%
%%%%%%%%%%%%%%%%%%%%%%%%%%%%%%%%%%%%%%%%%%%%%%%%%%%%%%%%%%%%%%%%%%%%%%%%%%%

\section{The direct spectral problem}

We will omit proofs which are nearly repetitions of the regular
case (for details see \cite{moi-these}).
\subsection{Solutions Properties}

In this section, $V$ is any $2\times 2$ matrix with $\lc$
coefficients. A fundamental system of solutions for
\eqref{akns-AKNS+V} when $V=0$ is given by
$$R(x,\lambda)=\frac{1}{\lambda^a}\Bigg[\begin{array}{c}
  \ja{a-1}{x} \\
  -\ja{a}{x}\\
\end{array}\Bigg],\quad\displaystyle S(x,\lambda)=\lambda^a\Bigg[\begin{array}{c}
  -\na{a-1}{x} \\
  \na{a}{x}\\
\end{array}\Bigg], $$
where $j_a$ and $\eta_a$ are spherical Bessel functions (see
section \ref{annexe-bessel}). These functions are called
fundamental since their wronskian is equal to $1$. From their
behavior near $x=0$, $R(x,\lambda)$ is called the regular
solution, it is analytic on $[0,1]\times\C$; $S(x,\lambda)$ is
called the singular solution, it is analytic on $(0,1]\times\C$.

Following Blancarte, Grébert and Weder \cite{bgw}, we construct
solutions for \eqref{akns-AKNS+V} by a Picard's iteration method
from $R$ and $S$.

Let $\sr$ and $\tilde{\ss}$ be defined by
$$\sr(x,\lambda,V)=\sum_{k\geq 0}R_k(x,\lambda,V),\quad\tilde{\ss}(x,\lambda,V)=\sum_{k\geq 0}S_k(x,\lambda,V)$$
with
\begin{eqnarray}
    \label{akns-relation_recurrence_reg}\left\{\begin{array}{c}
   R_0(x,\lambda,V)=R(x,\lambda), \\
      R_{k+1}(x,\lambda,V)=\displaystyle\int_0^x
    \G(x,t,\lambda)V(t)R_k(t,\lambda,V)dt,\quad k\in\N;
\end{array}\right.\\\label{akns-relation_recurrence_sing}
\left\{\begin{array}{c}
    S_0(x,\lambda,V)=S(x,\lambda), \\
    S_{k+1}(x,\lambda,V)=\displaystyle-\int_x^1
    \G(x,t,\lambda)V(t)S_k(t,\lambda,V)dt,\quad k\in\N.
\end{array}\right.
\end{eqnarray}
$\G$ is called Green function and is given by (see \cite{b})
\begin{equation}\label{akns-resolvante}
    \G(x,t,\lambda)=S(x,\lambda)R(t,\lambda)^\top-R(x,\lambda)S(t,\lambda)^\top.
\end{equation}

This construction is justified with the following
%%%%%%%%%%%%%%%%%%%%%%%%%%%%%%%%%%%%%%%%%%%%%%%%%%%%%%%%%%%%%%%%%%%%%%%%%%%
\begin{lem}\label{akns-lemme_base}
    Series defined by \eqref{akns-relation_recurrence_reg}, respectively
    by
    \eqref{akns-relation_recurrence_sing}, uniformly converge on
    bounded sets of $[0,1]\times\C\times\left(\lc\right)^4$, respectively
    of $(0,1]\times\C\times\left(\lc\right)^4$, towards solutions of \eqref{akns-AKNS+V}.
    Moreover, they satisfy the integral equations
    \begin{eqnarray*}
        \sr(x,\lambda,V)=R(x,\lambda)+\int_0^x
            \G(x,t,\lambda)V(t)\sr(t,\lambda,V)dt,\\
        \tilde{\ss}(x,\lambda,V)=S(x,\lambda)-\int_x^1
            \G(x,t,\lambda)V(t)\tilde{\ss}(t,\lambda,V)dt,
    \end{eqnarray*}
    and the estimates
    \begin{eqnarray*}
      \left|\sr(x,\lambda,V)\right| &\leq& C e^{|\im
            \lambda|x}\left(\frac{x}{1+|\lambda|x}\right)^a , \\
      \left|\tilde{\ss}(x,\lambda,V)\right|  &\leq& C e^{|\im
            \lambda|(1-x)}\left(\frac{1+|\lambda|x}{x}\right)^a,
    \end{eqnarray*}
with $C$ uniform on bounded sets of $\left(\lc\right)^4$.
\end{lem}
%%%%%%%%%%%%%%%%%%%%%%%%%%%%%%%%%%%%%%%%%%%%%%%%%%%%%%%%%%%%%%%%%%%%%%%%%%%
\noindent\textit{Proof.} We give it for $\sr$, it is similar for
$\tilde{\ss}$. Estimate
 \eqref{annexe-bessel_est_ja} for Bessel functions gives
\begin{equation}\label{akns-estimation-R}
    |R(x,\lambda)|\leq C e^{|\im\lambda|x}
    \left(\frac{x}{1+|\lambda|x}\right)^a.
\end{equation}
Iterative relation \eqref{akns-relation_recurrence_reg} leads to
\begin{equation}\label{akns-exp_R1}
    R_1(x,\lambda,V)=\int_0^x \G(x,t,\lambda)V(t)R(t,\lambda)dt,
\end{equation}
which, combining \eqref{akns-estimation-R} and the Green function
estimates \eqref{annexe-resolvante_est_tx}, is bounded by
\begin{equation*}
    \left|R_1(x,\lambda,V)\right|\leq C^2 e^{|\im\lambda|x}\left(\frac{x}{1+|\lambda|x}\right)^a \int_0^x |V(t)|dt,
\end{equation*}
By successive iterations and recurrence, for all positive integer
$n$, we get
\begin{equation*}
    \left|R_n(x,\lambda,V)\right|\leq \frac{C^{n+1}}{n!} e^{|\im\lambda|x}\left(\frac{x}{1+|\lambda|x}\right)^a
     \left(\int_0^{x} |V(t)|dt\right)^n.
\end{equation*}
This proves uniform convergence on bounded sets of
$[0,1]\times\C\times\left(\lc\right)^4$ for $\sr$ and the
estimate. Integral equation follows from
\eqref{akns-relation_recurrence_reg}. \qed

This uniform convergence gives us the following
\begin{prop}[Analyticity of solutions]$\quad$
\begin{enumerate}[(a)]
    \item  For all $x\in[0,1]$, $\sr(x,\lambda,V)$
        is analytic on $\C\times\left(\lc\right)^4$. Moreover, it is
        real valued on $\R\times\left(\lr\right)^4$.
    \item The map $\sr: (\lambda,V)\mapsto \sr(\cdot,\lambda,V)$ is analytic from $\C\times\left(\lc\right)^4$ to
    $H^1([0,1],\C^2)$.
    \item For all $x\in(0,1]$, $\tilde{\ss}(x,\lambda,V)$ is analytic on $\C\times\left(\lc\right)^4$ and real valued on
    $\R\times\left(\lr\right)^4$.
\end{enumerate}
\end{prop}

Let $\wr(\lambda,V)$ be the wronskian of $\sr$ and $\tilde{\ss}$,
defined by:
$$\wr(\lambda,V):=\wr\left(\sr(x,\lambda,V),\tilde{\ss}(x,\lambda,V)\right)=\det\left(\sr(x,\lambda,V),\tilde{\ss}(x,\lambda,V)\right).$$
Recall that $\wr(\lambda,V)$ is independent of $x$. We follow the
construction of a similar solution by Guillot and Ralston in
\cite{gr}: $\wr(\lambda,V)$ is not equal to $1$. However, as we
will see further, for $|\lambda|$ large enough, $\wr$ doesn't
vanishes (see Theorem \ref{akns-theoreme_asymptotique_L2_sing}).
Thus we may define the so-called singular solution by
$$\ss(x,\lambda,V)=\frac{\tilde{\ss}(x,\lambda,V)}{\wr(\lambda,V)},\quad x\in(0,1].$$
%%%%%%%%%%%%%%%%%%%%%%%%%%%%%%%%%%%%%%%%%%%%%%%%%%%%%%%%%%%%%%%%%%%%%%%%%%%
%\begin{rem}
%\noindent\textit{Remark.} %Here, we made a slight abuse of
%notation. Indeed, $\lambda\mapsto \wr(\lambda,V)$ may have zeros.
%But, first: since it is a analytic function, its zeros are
%isolated, so around each of them we just switch $\ss$ to another
%function independent from $\sr$. It is possible since the
%restriction of the ODE \eqref{akns-AKNS+V} to any subset
%$[\delta,1]$ with $\delta>0$ has a two dimensional space of
%solutions (see for instance \cite{gg}). Second, as we will see
%further, for $|\lambda|$ large enough, $\wr$ doesn't vanishes (see
%Theorem \ref{akns-theoreme_asymptotique_L2_sing}). The choice of
%this singular solution is justified by its asymptotic utility when
%$|\lambda|\rightarrow\infty$.\\
%Here, we made a slight abuse of notation. Indeed, $\lambda\mapsto
%\wr(\lambda,V)$ may have zeros. But, first: since it is a analytic
%function, its zeros are isolated, so around each of them we just
%switch $\ss$ to another function independent from $\sr$. It is
%possible since the restriction of the ODE \eqref{akns-AKNS+V} to
%any subset $[\delta,1]$ with $\delta>0$ has a two dimensional
%space of solutions (see for instance \cite{gg}). Second, . The
%choice of this singular solution is justified by its asymptotic
%utility when
%$|\lambda|\rightarrow\infty$.\\

%\end{rem}

Regularity of $\sr$ leads to existence of derivatives, obtained
following \cite{ist}:

\begin{prop}
For all $v\in\left(\lc\right)^4$, we have
\begin{eqnarray}
    \label{akns-differentielle_R_V}\left[d_{V} \sr(x,\lambda,V)\right](v)=\int_0^x
    \tilde{\G}(x,t,\lambda,V) v(t)\sr(t,\lambda,V) dt,\\
    \label{akns-gradient_R_lambda}\frac{\partial \sr}{\partial
    \lambda}(x,\lambda,V)=-\left[d_{V}
    \sr(x,\lambda,V)\right](\id),
\end{eqnarray}
where
$$\tilde{\G}(x,t,\lambda,V)=\ss(x,\lambda,V) \sr(t,\lambda,V)^\top-
\sr(x,\lambda,V)\ss(t,\lambda,V)^{\top}.$$
\end{prop}

\begin{nota}\label{akns-notation_Y1Z1ab}
For simplicity, we name the components of
solutions by
$$\sr(x,\lambda,p,q)=\Bigg[\begin{array}{c}
  Y_1(x,\lambda,p,q) \\
  Z_1(x,\lambda,p,q) \\
\end{array}\Bigg],\,\ss(x,\lambda,p,q)=\Bigg[\begin{array}{c}
  Y_2(x,\lambda,p,q) \\
  Z_2(x,\lambda,p,q) \\
\end{array}\Bigg]$$
and we introduce the following quantities
\begin{eqnarray*}
      a(x,\lambda,p,q)&= -\left[Y_1(x,\lambda,p,q)Z_2(x,\lambda,p,q)+Z_1(x,\lambda,p,q)Y_2(x,\lambda,p,q)\right], \\
      b(x,\lambda,p,q)&=\phantom{-}
      \left[Y_1(x,\lambda,p,q)Y_2(x,\lambda,p,q)-Z_1(x,\lambda,p,q)Z_2(x,\lambda,p,q)\right].
\end{eqnarray*}
\end{nota}

Now precise derivative expressions for AKNS potentials defined by
\eqref{akns-V_akns}. First, we define $\lc$-gradients for multiple
variable functions.

\begin{defn}
Let $H$ be an Hilbert space. For a continuously differentiable
complex valued map $f:(p,q)\mapsto f(p,q)$, the $\lc$-gradient
with respect to $(p,q)$ is the vector valued function
$$\displaystyle\nabla_{p,q}f=\left(\frac{\partial f}{\partial
p},\frac{\partial f}{\partial q}\right)$$ where
$\displaystyle\frac{\partial f}{\partial p}$, resp.
$\displaystyle\frac{\partial f}{\partial q}$ is the Riesz
representant of the partial differential $D_{p}f$, resp. $D_{q}f$
defined by
$$d_{p,q}f (v_1,v_2)= D_{p}f (v_1)+D_{q}f(v_2),
\quad (v_1,v_2)\in H\times H.$$
\end{defn}
\noindent\textit{Remark.} If $f$ is valued in $\C^n$, this
notation is understood component by component.\\

\begin{coro}[AKNS Gradients] For all $(p,q)\in\lc$, we have
\begin{eqnarray}
    \fl\qquad\left[\frac{\partial\sr}{\partial
    p}(x,\lambda,p,q)\right](t)=\indic_{[0,x]}(t)\Big[\ss(x,\lambda,p,q)\left[2Y_1(t,\lambda,p,q)Z_1(t,\lambda,p,q)\right]\nonumber\\
    +\sr(x,\lambda,p,q)a(t,\lambda,p,q)\Big],\label{akns-gradient_R_p}\\
    \fl\qquad \left[\frac{\partial\sr}{\partial
    q}(x,\lambda,p,q)\right](t)=\indic_{[0,x]}(t)\Big[\ss(x,\lambda,p,q)\left[Z_1(t,\lambda,p,q)^2-Y_1(t,\lambda,p,q)^2\right]\nonumber\\
    +\sr(x,\lambda,p,q)b(t,\lambda,p,q)\Big],\label{akns-gradient_R_q}\\
    \fl\qquad
   \left[\frac{\partial\sr}{\partial
    \lambda}(x,\lambda,p,q)\right]=\int_0^x\Big[-\ss(x,\lambda,p,q)\left[Y_1(t,\lambda,p,q)^2+Z_1(t,\lambda,p,q)^2\right]\nonumber\\
    +\sr(x,\lambda,p,q)\left[Y_1(t,\lambda,p,q)Y_2(t,\lambda,p,q)+Z_1(t,\lambda,p,q)Z_2(t,\lambda,p,q)\right]\Big]dt. \label{akns-gradient_R_lambda2}
\end{eqnarray}
\end{coro}
%   Propriété des valeurs propres

\subsection{Spectra}

%Considérons maintenant des potentiels réels.

Condition at $x=0$ selects a solution collinear to $\sr$,
condition at $x=1$ reduces spectrum to an eigenvalues-sequence. To
this end, we set:
\begin{nota}
Let $D(\lambda,V)$ be defined by:
\begin{equation}\label{akns-def-D}
    D(\lambda,V)=\sr(1,\lambda,V)\cdot u_\beta.
\end{equation}
Moreover, for all $u=(a,b)\in\C^2$, we define $u^\perp$ by
\begin{equation}\label{akns-def-perp}
    (a,b)^\perp=(b,-a).
\end{equation}
%En particulier, si $u$ est réel, $u$ et $u^\perp$ sont
%orthogonaux.
\end{nota}%
%Nous pouvons dès lors obtenir une nouvelle formulation pour les
%valeurs propres du problème:
\begin{prop}
    $D$ is analytic in $\lambda$ and $V$. The roots of $\lambda \mapsto D(\lambda,V)$ are exactly the eigenvalues
    for \eqref{akns-AKNS+V}-\eqref{akns-AKNS-bord}. Moreover, if $V$ is real-valued, they are all simple.
\end{prop}
\noindent\textit{Proof.} Analyticity of $D$ comes from $\sr$.
Since $\left\{\sr,\ss\right\}$ is a basis for the solutions of
\eqref{akns-AKNS+V}, the identification between eigenvalues and
roots of $\lambda\mapsto D(\lambda,V)$ follows.\\
Now suppose $V$ is real-valued et let $\lambda_0$ be an eigenvalue
of the problem. Simplicity lies on
\begin{equation}\label{akns-eq_spectre_simple}
    \n{\sr(\cdot,\lambda_0,V)}{\lr}^2=-(\sr(1,\lambda_0,V)\cdot
{u_\beta}^\perp)\frac{\partial D}{\partial \lambda}(\lambda_0,V).
\end{equation}
Indeed, from \eqref{akns-differentielle_R_V} and
\eqref{akns-gradient_R_lambda} we have
$$\frac{\partial D}{\partial \lambda}(\lambda_0,V)=-\left( \ss(1,\lambda_0,V)\cdot u_\beta\right)
\n{\sr(\cdot,\lambda_0,V)}{\lr}^2.$$ Then, rewriting the wronskian
of $\sr(1,\lambda_0,V)$ and $\ss(1,\lambda_0,V)$ in the
orthonormal basis $\left\{u_\beta,{u_\beta}^\perp\right\}$, we
obtain $\left( \sr(1,\lambda_0,V)\cdot
{u_\beta}^\perp\right)\left(\ss(1,\lambda_0,V)\cdot
u_\beta\right)=1$. \qed

From now, $V$ is defined by \eqref{akns-V_akns}, corresponding to
an AKNS operator.

%   Estimations H1

\subsection{$\hunc$-estimates}

In order to obtain accurate asymptotics, we add some regularity on
potentials. We use this roundabout method not because of the
singularity $a/x$ in the equation, but because of the AKNS
operator itself. Indeed, contrary to the Schrödinger operator,
there is no explicit decreasing for the Green function $\G$ with
respect to $\lambda$; so we have to force it allowing some
derivation. For the regular case ($a=0$), see for instance
\cite{gg}.

%%%%%%%%%%%%%%%%%%%%%%%%%%%%%%%%%%%%%%%%%%%%%%%%%%%%%%%%%%%%%%%%%%%%%%%%%%
\begin{thm}\label{akns-thm_estimation_h1_reg}
    For $(p,q)\in\left(\hunc\right)^2$, we have
    \begin{eqnarray}
    \fl\quad \Big|\sr(x,\lambda,p,q)-R(x,\lambda)\Big|\leq C\n{V}{\hunc}\left[\frac{x}{1+|\lambda x|}\right]^{a+1}
    \ln{[2+|\lambda x|]}e^{|\im
    \lambda|x+C\n{V}{2}},\label{akns-estimation_h1_reg}
\end{eqnarray}
uniformly on $[0,1]\times\C\times\left(\hunc\times\hunc\right)$,
where $\n{V}{\hunc}^2=\n{p}{\hunc}^2+\n{q}{\hunc}^2.$
\end{thm}
%%%%%%%%%%%%%%%%%%%%%%%%%%%%%%%%%%%%%%%%%%%%%%%%%%%%%%%%%%%%%%%%%%%%%%%%%%
\noindent\textit{Proof.} From relation
\eqref{akns-relation_recurrence_reg} at $k=1$ and
\eqref{akns-resolvante}, %the relation $$
%R_1(x,\lambda,p,q)=\int_0^x \G(x,t,\lambda)V(t)R_0(t,\lambda)d t$$
%becomes :
we have
\begin{eqnarray*}
 \fl R_1(x,\lambda,p,q)&=S_0(x,\lambda)\int_0^x R_0(t,\lambda)^\top V(t)R_0(t,\lambda)d t\\
   &\phantom{=\,}-R_0(x,\lambda)\int_0^x S_0(t,\lambda)^\top V(t)R_0(t,\lambda)d
  t\\
&=S_0(x,\lambda)\int_0^x \left[q(t)\left(
R_0^2(t,\lambda)^2-R_0^1(t,\lambda)^2\right)%\\
%&\qquad\qquad\qquad
+2p(t)R_0^1(t,\lambda)R_0^2(t,\lambda)\right] d t\\
&\phantom{=} -R_0(x,\lambda)\int_0^x \left[q(t)\left(
S_0^2(t,\lambda)R_0^2(t,\lambda)-S_0^1(t,\lambda)R_0^1(t,\lambda)\right)\right.\\
& \phantom{= -R_0(x,\lambda)\int_0^x \,}
\left.+p(t)\left(S_0^1(t,\lambda)R_0^2(t,\lambda)+S_0^2(t,\lambda)R_0^1(t,\lambda)\right)\right]
d
  t.
\end{eqnarray*}
We can write
$R_1(x,\lambda,p,q)=\lambda^{-a}\left[X(q)+Y(p)\right],$ where
\begin{eqnarray*} \fl X(q)=\Bigg[\begin{array}{c}
  \displaystyle-\na{a-1}{x} \\
 \displaystyle \na{a}{x}\\
\end{array}\Bigg]\int_0^x \Big\{[\ja{a}{t}]^2-[\ja{a-1}{t}]^2\Big\}q(t) d
t\\
+\Bigg[\begin{array}{c}
\displaystyle\ja{a-1}{x}\\
\displaystyle-\ja{a}{x}\\
\end{array}\Bigg]\int_0^x \Big[\na{a}{t}\ja{a}{t}-\na{a-1}{t}\ja{a-1}{t}\Big]q(t)
dt,\\
\fl Y(p)=\Bigg[\begin{array}{c}
  \displaystyle \na{a-1}{x} \\
 \displaystyle -\na{a}{x}\\
\end{array}\Bigg]\int_0^x \Big[2\ja{a-1}{t}\ja{a}{t}\Big]p(t) d
t\\
-\Bigg[\begin{array}{c}
\displaystyle\ja{a-1}{x}\\
\displaystyle-\ja{a}{x}\\
\end{array}\Bigg]\int_0^x \Big[\na{a-1}{t}\ja{a}{t}+\na{a}{t}\ja{a-1}{t}\Big]p(t) dt.
\end{eqnarray*}
%%%%%%%%%%%%%%%%%%%%%%%%%%%%%%%%%%%%%%%%%%%%%%%%%%%%%%%%%%%%%%%%%%%%%%%%%%%%%%%%%%%%%%%%%%%%%%
\begin{description}

%   Majoration de X

    \item[Estimation for $X(q)$:]$\,$\\
    Integrating by parts, we get
    \begin{eqnarray*} \fl X(q)= \frac{1}{\lambda}\Bigg[\begin{array}{c}
        0 \\ -\ja{a}{x} \\ \end{array}\Bigg]    q(x)\\
        +\frac{1}{\lambda}\int_0^x\Bigg[\begin{array}{c}
        -\na{a-1}{x}\ja{a-1}{t}+\ja{a-1}{x}\na{a-1}{t} \\
        \na{a}{x}\ja{a-1}{t}-\ja{a}{x}\na{a-1}{t}\\
    \end{array}\Bigg]\ja{a}{t}q'(t)dt.
    \end{eqnarray*}
    Estimates \eqref{annexe-bessel_est_ja},\eqref{annexe-resolvante_est_tx} and Sobolev inequality %(cf. \cite{br})
    $\n{q}{\infty}\leq  C\n{q}{\hunc}$
    %th VIII.7 p.129: Il existe une constante $C$ telle que $$\n{f}{\infty}\leq C \n{f}{W^{1,p}}\quad\forall f\in W^{1,p}([0,1]),\quad
    %\forall p\in [1,\infty].$$}
    give \begin{equation}\label{akns-estimation_Xq}
        |X(q)|\leq\frac{C }{|\lambda|}\left(\frac{|\lambda
        x|}{1+|\lambda x|}\right)^{a+1}e^{|\im
        \lambda|x}\n{q}{\hunc}.
    \end{equation}

%   Majoration de Y

    \item[Estimation for $Y(p)$:]$\,$\\
    With notations from lemmas \ref{annexe-lemme_primitive_K} and \ref{annexe-lemme_primitive_L},
    integration by parts gives :
        \begin{eqnarray*}\fl Y(p)=\Bigg[\begin{array}{c}
            \displaystyle \na{a-1}{x} \\
            \displaystyle -\na{a}{x}\\
            \end{array}\Bigg]\left(\left[\frac{1}{\lambda}F_1(\lambda
            t)p(t)\right]_0^x-\frac{1}{\lambda}\int_0^x F_1(\lambda t)p'(t) dt\right)\\
            -\Bigg[\begin{array}{c}\displaystyle\ja{a-1}{x}\\
            \displaystyle-\ja{a}{x}\\
            \end{array}\Bigg]
            \left(\left[\frac{1}{\lambda}F_2(\lambda
            t)p(t)\right]_0^x-\frac{1}{\lambda}\int_0^x F_2(\lambda t)p'(t)
            dt\right).
        \end{eqnarray*}
    %Nous pouvons maintenant déterminer les estimations sur $Y(p)$.
%\begin{description}
    \item[$\quad$ When $|\lambda x|\leq$1.]%

    Estimations  \ref{annexe-bessel_est_ja}, \ref{annexe-bessel_est_na} and part \textit{(i)}
    from lemmas \ref{annexe-lemme_primitive_K} and
    \ref{annexe-lemme_primitive_L} lead to
    \begin{eqnarray*} |Y(p)|\leq&\frac{C }{|\lambda|}\left(\frac{|\lambda x|}{1+|\lambda
        x|}\right)^{a+1}\n{p}{\hunc}e^{|\im \lambda|x}.
    \end{eqnarray*}

    \item[$\quad$ When $|\lambda x|\geq$1.]%

    Now, we only consider $Y(p)$ second component, the proof is similar for the first one. Terms to estimate
    contain :
    $$g(x,t):=\na{a}{x}F_1(\lambda t)-\ja{a}{x}F_2(\lambda
    t),\quad0\leq t\leq x.$$

%    \begin{description}

        \item[$\qquad$ If $|\lambda t|\leq$1.]
        As for $|\lambda x|\leq 1$, we get
        $\displaystyle
        |g(x,t)|\leq 2C\left(\frac{|\lambda x|}{1+|\lambda
            x|}\right)^{a+1}e^{|\im \lambda|x}.
        $
        \item[$\qquad$ If $|\lambda t|\geq$1.]
        Using points  \textit{(ii)} from lemmas \ref{annexe-lemme_primitive_K}
        and \ref{annexe-lemme_primitive_L}, expressions \eqref{annexe-bessel-ja-sinus}
        and
        \eqref{annexe-bessel-na-sinus}, it follows :
        \begin{eqnarray*}
        \fl\qquad    g(x,t)=\eta_a(\lambda x)r_a(\lambda t)\\
            -a \left[\cos{\left(\lambda x-\frac{a\pi}{2}\right)}\mathrm{ci}{(2\lambda
            t)}+\sin{\left(\lambda x-\frac{a\pi}{2}\right)}\mathrm{Si}{(2\lambda t)}
            \right]P_a(\lambda x)\\
            +a \left[\sin{\left(\lambda x-\frac{a\pi}{2}\right)}\mathrm{ci}{(2\lambda
            t)}-\cos{\left(\lambda x-\frac{a\pi}{2}\right)}\mathrm{Si}{(2\lambda t)}
            \right]I_a(\lambda x)\\
            +\left(P_a(\lambda x)p_a(\lambda t)-I_a(\lambda x)q_a(\lambda t) \right)\cos{\left[\lambda (x-2t)-\frac{a\pi}{2}\right]}\\
            -\left(P_a(\lambda x)q_a(\lambda t)+I_a(\lambda x)p_a(\lambda t) \right)\sin{\left[\lambda
            (x-2t)-\frac{a\pi}{2}\right]}.
        \end{eqnarray*}
        (To lighten, the polynomial variable $X$ is replaced by
        $1/X$.) First term is bounded by $C
        e^{|\im\lambda|x}$ thanks to \eqref{annexe-bessel_est_na}.
        The last two terms are uniformly bounded by $C e^{|\im \lambda|(x-2t)}$ on the considered area.
        Now remains the following expression
        $$h(x,t):=\cos{\left(\lambda x-\frac{a\pi}{2}\right)}\mathrm{ci}{(2\lambda
            t)}+\sin{\left(\lambda x-\frac{a\pi}{2}\right)}\mathrm{Si}{(2\lambda t)}.$$
        According to \cite{wolf-eq1} and \cite{wolf-eq2}, we have
        \begin{eqnarray*}
        \fl\qquad  \mathrm{ci}(z) = -\gamma-\frac{\log(z^2)}{2}+\frac{\sin{z}}{z}\left(1+\mathcal{O}_1\left(\frac{1}{z^2}\right)\right)%\\
          %&& \qquad
          -\frac{\cos{z}}{z^2}\left(1+\mathcal{O}_2\left(\frac{1}{z^2}\right)\right),\\
        \fl\qquad  \mathrm{Si}(z) =
          \frac{\pi\sqrt{z^2}}{2z}-\frac{\cos{z}}{z}\left(1+\mathcal{O}_1\left(\frac{1}{z^2}\right)\right)%\\
          %&& \qquad
          -\frac{\sin{z}}{z^2}\left(1+\mathcal{O}_2\left(\frac{1}{z^2}\right)\right).
        \end{eqnarray*}
        Thus, we get
        \begin{eqnarray*}
        \fl\qquad   h(x,t)=-\left[\gamma+\frac{\log{(2\lambda t)^2}}{2}\right]
                    \cos{\left(\lambda x-\frac{a\pi}{2}\right)}
                    %\\
            +\frac{\pi\sqrt{(2\lambda t)^2}}{2\lambda t}
                    \sin{\left(\lambda x-\frac{a\pi}{2}\right)}\\
            -\frac{1}{2\lambda t}\left(1+\mathcal{O}_1\left(\frac{1}{(2\lambda
            t)^2}\right)\right)\sin{\left[\lambda(x-2t)-\frac{a\pi}{2}\right]}\\
            -\frac{1}{(2\lambda t)^2}\left(1+\mathcal{O}_2\left(\frac{1}{(2\lambda
            t)^2}\right)\right)
            \cos{\left[\lambda (x-2t)-\frac{a\pi}{2}\right]}.
        \end{eqnarray*}
    The last three terms are also uniformly controlled by $C e^{|\im \lambda|(x-2t)}$;
    the first one is bounded by $C \ln{|\lambda t|} e^{|\im \lambda|x}$.
%    \end{description}
%\end{description}
Combining the above estimates, we obtain the following uniform
estimate
\begin{equation}\label{akns-estimation_Yp}
    |Y(p)|\leq\frac{C}{|\lambda|}\left(\frac{|\lambda x|}{1+|\lambda
    x|}\right)^{a+1}\ln{[2+|\lambda|x]}\n{p}{\hunc}e^{|\im
    \lambda|x}.
\end{equation}
\end{description}

Relations \eqref{akns-estimation_Xq}-\eqref{akns-estimation_Yp}
and the concavity rule
$$\forall (x,y)\in\R^2, \quad \frac{|x|+|y|}{2}\leq \sqrt{\frac{|x|^2+|y|^2}{2}},$$
imply
\begin{equation}\label{akns-estimation_R1_H1}
        \fl\qquad |R_1(x,\lambda,p,q)|\leq\frac{C
    }{|\lambda|^{a+1}}\left(\frac{|\lambda x|}{1+|\lambda
    x|}\right)^{a+1}\ln{[2+|\lambda|x]}\n{V}{\hunc}e^{|\im
    \lambda|x}.
\end{equation}
From this estimate, as in the proof of lemma
\ref{akns-lemme_base}, we deduce estimate
\eqref{akns-estimation_h1_reg}. \qed
%%%%%%%%%%%%%%%%%%%%%%%%%%%%%%%%%%%%%%%%%%%%%%%%%%%%%%%%%%%%%%%%%%%%%%%%%%%%%%%%%%%

\noindent\textit{Remark.}\label{remarque-dirac} A similar
computation for the Dirac operator is not easy, even if $a=0$.
Indeed, when we compute the term $R_1$, we do not only get a term,
loosely speaking, in $\O{1/\lambda}$ but also in $\O{1}$. And when
iterating this, we get at each time a new term $\O{1}$ and
$\O{1/\lambda}$. A way through this problem is given in
\cite{moi-these} using the latter gauge transformation to deduce,
for any $a\in\N$, some partial results from AKNS to Dirac
operator: spectrum, asymptotic expansion for eigenvalues and
eigenvectors, Borg-Levinson theorem type\ldots

%   Estimations H1

\subsection{$\lc$-Estimates}

To transform $\hunc$-estimates into $\lc$-estimates, we need an
auxiliary lemma (for the regular case, see \cite{ag}, \cite{mis}
and \cite{gkp}).
%%%%%%%%%%%%%%%%%%%%%%%%%%%%%%%%%%%%%%%%%%%%%%%%%%%%%%%%%%%%%%%%%%%%%%%%%%
\begin{lem}\label{akns-estimation_L2} Let $V_0\in\lc\times\lc$, $r_0\geq 0$, $\varepsilon\geq 0$ and let $V_\epsilon\in\hunc\times\hunc$
such that $\n{V_0-V_\varepsilon}{2}<\varepsilon$. Then, for all
$V\in\lc\times\lc$ such that $\n{V-V_0}{2}< r_0$ and for all
$(x,\lambda)\in[0,1]\times \C^\ast$, we have
\begin{eqnarray}
\fl\qquad    |\sr(x,\lambda,p,q)-R(x,\lambda)|\leq C&
\left(r_0+\varepsilon+\frac{\ln|\lambda|}{|\lambda|}\|V_\varepsilon\|_{\hunc}\right)\nonumber\\
    &\qquad\times
    \left(\frac{x}{1+|\lambda x|}\right)^a e^{|\im \lambda|x+C\n{V}{2}}.\label{akns-estimation-reguliere}
\end{eqnarray}
\end{lem}
%%%%%%%%%%%%%%%%%%%%%%%%%%%%%%%%%%%%%%%%%%%%%%%%%%%%%%%%%%%%%%%%%%%%%%%%%%
\noindent\textit{Proof.} Since $V_\epsilon\in\hunc\times\hunc$,
estimate \eqref{akns-estimation_R1_H1} obtained during the proof
of Theorem \ref{akns-thm_estimation_h1_reg} becomes
\begin{equation*}
\fl\qquad    |R_1(x,\lambda,V_\varepsilon)|\leq\frac{C
    }{|\lambda|^{a+1}}\left(\frac{|\lambda x|}{1+|\lambda
    x|}\right)^{a+1}\ln{[2+|\lambda|x]}\n{V_\varepsilon}{\hunc}e^{|\im
    \lambda|x}.
\end{equation*}
Using $\n{V_0-V_\varepsilon}{2}<\varepsilon$ and $\n{V-V_0}{2}<
r_0$ in \eqref{akns-exp_R1}, estimations
\eqref{annexe-resolvante_est_tx} and \eqref{annexe-bessel_est_ja}
lead to
\begin{equation*}
\fl\qquad
|R_1(x,\lambda,p,q)-R_1(x,\lambda,V_\varepsilon)|\leq\frac{C
    }{|\lambda|^{a}}\left(\frac{|\lambda x|}{1+|\lambda
    x|}\right)^{a}(r_0+\varepsilon)e^{|\im
    \lambda|x}.
\end{equation*}
Combining these two inequalities, we get
\begin{equation*}
\fl\qquad    |R_1(x,\lambda,p,q)|\leq C\left(\frac{x}{1+|\lambda
    x|}\right)^{a}\left(r_0+\varepsilon+\frac{\ln{[2+|\lambda|x]}}{|\lambda|}\n{V_\varepsilon}{\hunc}\right)e^{|\im
    \lambda|x}.
\end{equation*}
Iterating this with \eqref{akns-relation_recurrence_reg}, we
deduce for every $n\in\N$
\begin{eqnarray*}
\fl\qquad    |R_{n+1}(x,\lambda,p,q)|\leq \frac{C^{n+1}}{n!}\left(r_0+\varepsilon+\frac{\ln{[2+|\lambda|x]}}{|\lambda|}\n{V_\varepsilon}{\hunc}\right)\\
    \times\left(\frac{x}{1+|\lambda
    x|}\right)^{a}e^{|\im
    \lambda|x}\left(\int_0^1 |V(t)|dt\right)^n.
\end{eqnarray*}
Then, summing up, estimation \eqref{akns-estimation-reguliere}
follows. \qed

We now deduce the following %
%%%%%%%%%%%%%%%%%%%%%%%%%%%%%%%%%%%%%%%%%%%%%%%%%%%%%%%%%%%%%%%%%%%%%%%%%%
\begin{prop}\label{akns-prop_L2_loc} Let $(p,q)\in\left(\lc\right)^2$, we have uniformly on $[0,1]$,
\begin{equation}\label{akns-estimation_L2_loc}
\sr(x,\lambda,p,q)=R(x,\lambda)+o\left[\left(\frac{x}{1+|\lambda
x|}\right)^a e^{|\im \lambda|x}\right], \quad
|\lambda|\rightarrow\infty.
\end{equation}
\end{prop}
%%%%%%%%%%%%%%%%%%%%%%%%%%%%%%%%%%%%%%%%%%%%%%%%%%%%%%%%%%%%%%%%%%%%%%%%%%
\noindent\textit{Proof.} From Lemma \ref{akns-estimation_L2} with
$r_0=0$, given $\delta>0$ there exists $\lambda_\delta>0$ such
that
$$
|\sr(x,\lambda,p,q)-R(x,\lambda)|\leq \delta
\left(\frac{x}{1+|\lambda x|}\right)^a e^{|\im
\lambda|x+C\n{V}{2}},$$ for all $\lambda$ such that
$|\lambda|>\lambda_\delta$. \qed

%   Localisation du spectre

\subsection{Spectrum localization}

%%%%%%%%%%%%%%%%%%%%%%%%%%%%%%%%%%%%%%%%%%%%%%%%%%%%%%%%%%%%%%%%%%%%%%%%%%
\begin{thm}[Counting Lemma]\label{akns-lemme_comptage}$\,$\\
Let $(p_0,q_0)\in\lc\times\lc$, there exist $\varepsilon>0$ and an
integer $N_0>0$ such that for all $(p,q)\in\lc\times\lc$ with
$\n{(p,q)-(p_0,q_0)}{\lc}<\varepsilon$, the following statements
hold:
\begin{itemize}
    \item [$\bullet$]For all $|n|>N_0$, $\lambda\mapsto D(\lambda,p,q)$ has exactly one root in
    $\left|\lambda-\left(n\pi+\frac{a\pi}{2}+\beta\right)\right|<\frac{\pi}{2}$,
    \item [$\bullet$]$\lambda\mapsto D(\lambda,p,q)$ has exactly $2N_0+1-a$
    root counted with multiplicity in
    $\left|\lambda-\left(\frac{a\pi}{2}+\beta\right)\right|<\left(N_0+\frac{1}{2}\right)\pi$,
    \item [$\bullet$]$\lambda\mapsto D(\lambda,p,q)$ has no root elsewhere.
\end{itemize}
\end{thm}
%%%%%%%%%%%%%%%%%%%%%%%%%%%%%%%%%%%%%%%%%%%%%%%%%%%%%%%%%%%%%%%%%%%%%%%%%%
\noindent\textit{Proof.} Let $\varepsilon>0$, from estimate
\eqref{akns-estimation-reguliere}  and using Lemma
\ref{akns-estimation_L2} notations, we have
\begin{equation*}%
|\sr(1,\lambda,p,q)-R(1,\lambda)|\leq\frac{C e^{|\im
\lambda|+C\n{V}{2}}}{|\lambda|^{a}}\left(2\varepsilon+\frac{\ln{|\lambda|}}{|\lambda|}\n{V_\varepsilon}{\hunc}\right).
\end{equation*}%
Bessel functions relation \eqref{annexe-bessel-ja-sinus} implies
the following uniform estimate on $|\lambda|>1$
\begin{equation}\label{akns-est-R11}
    R(1,\lambda)=\frac{1}{\lambda^a}\Bigg[\begin{array}{c}
  \cos{\left(\lambda-\frac{a\pi}{2}\right)} \\
  -\sin{\left(\lambda-\frac{a\pi}{2}\right)} \\
\end{array}\Bigg]+\mathcal{O}\left(\frac{e^{|\im
\lambda|}}{|\lambda|^{a+1}}\right),
\end{equation}
which leads, together with the previous one, to
\begin{eqnarray*}
\fl    \left|\lambda^a
D(\lambda,p,q)-\sin\left(\beta+\frac{a\pi}{2}-\lambda\right)\right|\leq
\left(2\varepsilon+\frac{\ln{|\lambda|}}{|\lambda|}\n{V_\varepsilon}{\hunc}+\frac{1}{|\lambda|}\right)
 C e^{C\n{V}{2}}e^{|\im \lambda|}.
\end{eqnarray*}
 Now introduce the circles:\begin{itemize}
    \item for $n\in\Z$, $\gamma_n$ is defined by
     $$\left|\lambda-\left(n\pi+\frac{a\pi}{2}+\beta\right)\right|=\frac{\pi}{2}.$$
    \item for $n\in\N$, $C_n$ is defined by
    $$\left|\lambda-\left(\frac{a\pi}{2}+\beta\right)\right|=\left(n+\frac{1}{2}\right)\pi.$$
 \end{itemize}
\noindent We choose $\varepsilon>0$ such that $C e^{C\n{V}{2}}
2\varepsilon<\frac{1}{8}$. Moreover, on each circle we have
$$|\lambda|>\left(N_0+\frac{1}{2}\right)\pi-\left|\frac{a\pi}{2}+\beta\right|%\quad \mathrm{(voir le schéma ci-dessus)}
$$
 and since the map $t\mapsto\frac{\ln{t}}{t}$ decreases on $]e,\infty[$, we can pick up
$N_0>0$ such that
$$C e^{C\n{V}{2}}\frac{\ln{|\lambda|}}{|\lambda|}\n{V_\varepsilon}{\hunc}<\frac{1}{8}.$$

Thus, we get the following
$$\left|\lambda^a
D(\lambda,p,q)-\sin\left(\beta+\frac{a\pi}{2}-\lambda\right)\right|<\frac{1}{4}\,e^{|\im
\lambda|}=\frac{1}{4}\,e^{|\im
\left(\lambda-\frac{a\pi}{2}-\beta\right)|}.$$

Using the following estimate for all $k\in\Z$ (see Lemma 2.1 in
\cite{ist})
$$e^{|\im z|}< 4 |\sin{z}|\quad \mathrm{for} \quad|z-k\pi|\geq \frac{\pi}{4},$$
on the sets $\gamma_n$ and $C_{N_0}$,
 with $z=\lambda-\frac{a\pi}{2}-\beta$, we obtain
$$\left|\lambda^a
D(\lambda,p,q)-\sin\left(\beta+\frac{a\pi}{2}-\lambda\right)\right|<
\left|\sin\left(\beta+\frac{a\pi}{2}-\lambda\right)\right|.$$ Now,
the use of the Rouché Theorem let us conclude that the analytical
functions $\lambda\mapsto\lambda^a D(\lambda,p,q)$ and
$\lambda\mapsto\sin\left(\beta+\frac{a\pi}{2}-\lambda\right)$ have
the same number of roots counted with multiplicity inside theses
circles. To show there is no other elsewhere, we just have to
consider an other circle $C_N$ with $N>N_0$ and apply again the
Rouché Theorem. \qed

Now, we can order eigenvalues: when $n>N_0$, $\lambda_{a,n}(p,q)$
is the eigenvalue surrounded by $\gamma_n$. Next, we order
lexicographically the $2N_0+1-a$ eigenvalues lying in $C_{N_0}$,
in other words, for $k=a-N_0,\dots,N_0-1$:
\begin{eqnarray*}
    \re{\lambda_{a,k}(p,q)}<\re{\lambda_{a,k+1}(p,q)}\\
    \mathrm{ or }\\
    \re{\lambda_{a,k}(p,q)}=\re{\lambda_{a,k+1}(p,q)} \quad\mathrm{
    and
    }\quad\im{\lambda_{a,k}(p,q)}\leq\im{\lambda_{a,k+1}(p,q)}.
\end{eqnarray*}
To continue the numbering, the eigenvalue included in
$\gamma_{-n}$, for $n>N_0$, must be $\lambda_{a,-n+a}$. To put it
directly, we say that for $n>N_0-a$, $\lambda_{a,-n}$ is the
eigenvalue surrounded by $\gamma_{-(n+a)}$.

The localization gives us the following locally uniform estimates
on $\lc\times\lc$
\begin{eqnarray}\label{akns-asymptotique_vp_L2_1}
    \lambda_{a,n}(p,q)=\left(n+\frac{a}{2}\right)\pi+\beta+\mathcal{O}(1),\quad
    n\rightarrow\infty,\quad |\O{1}|\leq \frac{\pi}{2},\\\label{akns-asymptotique_vp_L2_2}
    \lambda_{a,-n}(p,q)=-\left(n+\frac{a}{2}\right)\pi+\beta+\mathcal{O}(1),\quad
    n\rightarrow\infty,\quad |\O{1}|\leq \frac{\pi}{2}.
\end{eqnarray}
%%%%%%%%%%%%%%%%%%%%%%%%%%%%%%%%%%%%%%%%%%%%%%%%%%%%%%%%%%%%%%%%%%%%%%%%%
\begin{prop} Let $(p,q)\in\lc\times\lc$.
\begin{equation}\label{akns-est_lambda_petit o}
\lambda_{a,n}(p,q)=\left(n+\sgn(n)\frac{a}{2}\right)\pi+\beta+o(1)\quad,\quad
|n|\rightarrow+\infty.
\end{equation}
\end{prop}
\noindent\textit{Proof.} Relation \eqref{akns-estimation_L2_loc}
at $x=1$ and definition \eqref{akns-def-D} give
$$D(\lambda,p,q)=R(1,\lambda)\cdot u_\beta+o\left(\frac{e^{|\im\lambda|}}{|\lambda|^a}\right)$$
then, estimate \eqref{akns-est-R11} implies
\begin{equation*}%\label{akns-vp-bidon-1}
    D(\lambda,p,q)=\frac{1}{\lambda^a}\sin{\left(\beta+\frac{a\pi}{2}-\lambda\right)}+o\left(\frac{e^{|\im\lambda|}}{|\lambda|^a}\right).
\end{equation*}
According to the counting lemma, we have
\begin{equation*}%\label{akns-vp-bidon-2}
    \lambda_{a,n}(p,q)=\left(n+\sgn(n)\frac{a}{2}\right)\pi+\beta+\O{1},\quad
    |n|\rightarrow\infty,
\end{equation*}
knowing that
$
%\begin{equation}\label{akns-vp-bidon-3}
    |\O{1}|<\frac{\pi}{2}.
%\end{equation}
$ We evaluate $D(\lambda,p,q)$ at $\lambda=\lambda_{a,n}(p,q)$ and
use the above estimates to get
\begin{equation*}%\label{akns-vp-bidon-4}
    0=\frac{1}{{\lambda_{a,n}}^a}\sin{(\O{1})}+o\left(\frac{1}{|\lambda_{a,n}|^a}\right).
\end{equation*}
By identification, we found the result. \qed

\noindent\textit{Remarks} Theses results have to be compared with
those in the regular case ($a=0$):
\begin{itemize}
    \item Asymptotics of solutions and eigenvalues localization for $\lc$-potentiels
    are only locally uniform. This is due to the operator by
    itself and not to the singularity. In \cite{gk} is given a pair of potentials with
    identical $\lc$-norm whose eigenvalues numbering (localization) are different.
    \item A new phenomenon, relative to the numbering, is this loss of $a$
    eigenvalues lying near $0$. It may be seen as the analogue of the shift by $a/2$ in the
    eigenvalues asymptotics of the radial Schrödinger operator (see for instance \cite{gr} when $a=1$ and \cite{moi-SL} for general $a$).
\end{itemize}

%%%%%%%%%%%%%%%%%%%%%%%%%%%%%%%%%%%%%%%%%%%%%%%%%%%%%%%%%%%%%%%%%%%%%%%%%%%%%%%%%%%%%%%%%%%%%%%%%%%%%%%%%
%%%%%%%%%%%%%%%%%%%%%%%%%%%%%%%%%%%%%%%%%%%%%%%%%%%%%%%%%%%%%%%%%%%%%%%%%%%%%%%%%%%%%%%%%%%%%%%%%%%%%%%%%
\section{Spectral Data}
%%%%%%%%%%%%%%%%%%%%%%%%%%%%%%%%%%%%%%%%%%%%%%%%%%%%%%%%%%%%%%%%%%%%%%%%%%%%%%%%%%%%%%%%%%%%%%%%%%%%%%%%%
%%%%%%%%%%%%%%%%%%%%%%%%%%%%%%%%%%%%%%%%%%%%%%%%%%%%%%%%%%%%%%%%%%%%%%%%%%%%%%%%%%%%%%%%%%%%%%%%%%%%%%%%%
From this point, $(p,q)$ are real-valued. Thus,
$(\lambda_{a,n}(p,q))_{n\in\Z}$ is a strictly increasing sequence
of real numbers. We set some notations  :
\begin{nota} We define
$$\sr_n(t,p,q)=\sr(t,\lambda_{a,n}(p,q),p,q)\quad\mathrm{ and }\quad\ss_n(t,p,q)=\ss(t,\lambda_{a,n}(p,q),p,q).$$
Let $G_n(t,p,q)$ be the normed eigenvector with respect to
$\lambda_{a,n}(V)$ defined by
$$G_n(t,p,q)=\frac{\sr_n(t,p,q)}{\n{\sr_n(\cdot,p,q)}{2}}.$$
We also define
$$A_n(x,p,q)=\left(a_n(x,p,q),b_n(x,p,q)\right)$$
where $a_n(x,p,q)= a(x,\lambda_{a,n}(p,q),p,q)$ and $b_n(x,p,q)=
b(x,\lambda_{a,n}(p,q),p,q)$ ($a$ and $b$ are given on page
\pageref{akns-notation_Y1Z1ab}).
\end{nota}

\subsection{Regularity, derivatives}

Eigenvalues regularity and associated derivatives follows like in
\cite{gg} and \cite{ist} as pictured by the next proposition.
%%%%%%%%%%%%%%%%%%%%%%%%%%%%%%%%%%%%%%%%%%%%%%%%%%%%%%%%%%%%%%%%%%%%%%%%%%%
\begin{prop}
    For all $n\in\Z$, $(p,q)\mapsto \lambda_{a,n}(p,q)$ is a real-analytic map on $\lr\times\lr$. Its
    $\lr$-gradient is given by
    \begin{equation}\label{akns-expr_grad_lambda}
\fl \qquad
\nabla_{p,q}\lambda_{a,n}=\left(\frac{\partial\lambda_{a,n}}{\partial
    p},\frac{\partial\lambda_{a,n}}{\partial q} \right) \mathrm{ with }
    \left\{\begin{array}{l}  \displaystyle\frac{\partial\lambda_{a,n}}{\partial
    p}=2\, G_{n,1}(t,p,q)\,G_{n,2}(t,p,q), \\
    \displaystyle \frac{\partial\lambda_{a,n}}{\partial
    q}=G_{n,2}(t,p,q)^2-G_{n,1}(t,p,q)^2.
    \end{array}\right.
    \end{equation}
\end{prop}

%\noindent\textit{Proof.} First, relation
%\eqref{akns-eq_spectre_simple} gives
%$$\frac{\partial D}{\partial \lambda}(\lambda_{a,n}(p,q),p,q)\neq 0.$$
%Second, since $(p,q)\mapsto \sr(1,\lambda,p,q)$ is continuous, the
%counting lemma implies that $(p,q)\mapsto\lambda_{a,n}(p,q)$ is
%continuous. Thus, the implicit function theorem applied to $D$,
%real-analytic function on $\R\times\lr\times\lr$, gives the
%real-analyticity of $(p,q)\mapsto\lambda_{a,n}(p,q)$.\\
%To obtain its gradient, just notice that
%\begin{eqnarray*}
%    0&=\frac{\partial}{\partial
%    p}\left(D(\lambda_{a,n}(p,q),p,q)\right)\\
%    &=\frac{\partial \lambda_{a,n}}{\partial p}\left(\frac{\partial \sr}{\partial \lambda}(1,\lambda_{a,n},p,q)\cdot
%    u_\beta\right)+\frac{\partial \sr}{\partial
%    p}(1,\lambda_{a,n},p,q)\cdot u_\beta.
%\end{eqnarray*}
%Then, using relations \eqref{akns-gradient_R_p} and
%\eqref{akns-gradient_R_lambda2} , we get
%\begin{eqnarray*}
%\fl\qquad    0=-\frac{\partial \lambda_{a,n}}{\partial
%p}&\left(\ss(1,\lambda_{a,n},p,q)\cdot
%    u_\beta \right)\int_0^1 \left[Y_1(x,\lambda_{a,n},p,q)^2+Z_1(x,\lambda_{a,n},p,q)^2\right]dt\\
%    &+\left(\ss(1,\lambda_{a,n},p,q)\cdot
%    u_\beta \right) 2Y_1(t,\lambda_{a,n},p,q)Z_1(t,\lambda_{a,n},p,q).
%\end{eqnarray*}
%A similar computation gives the $q$-gradient and the result. \qed

%\subsection{Additional data}

Like in \cite{ist}, or simply following \cite{gg}, we need more
information to recover a complete parametrization of
$\left(\lr\right)^2$. Boundary condition at $x=1$ defining each
eigenvalue is an orthogonality relation following one direction.
It sounds reasonable  that the knowledge of a similar data in a
complementary (here orthogonal) direction is enough.
\begin{defn}
For all $n\in\Z$, we call normalization constants the quantities
\begin{equation}\label{akns-def_kappan}
    \kappa_{a,n}(p,q)=\sr_n(1,p,q)\cdot {u_\beta}^\perp.
\end{equation}
\end{defn}
Following \cite{gg}, we get :
\begin{prop}
    For all $n\in\Z$, $(p,q)\mapsto \kappa_{a,n}(p,q)$ is a real-analytic map on $\lr\times\lr$.
    Its $\lr$-gradient is given by
    \begin{equation}\label{akns-expr_grad_kappa}
    \frac{\nabla_{p,q}\kappa_{a,n}}{\kappa_{a,n}}=A_n(x,p,q)+\Big\langle
    \sr_n(\cdot,p,q),\ss_n(\cdot,p,q)\Big\rangle
    \nabla_{p,q}\lambda_{a,n}(p,q).\end{equation}
\end{prop}

Now, precise the behavior of theses normalization constants.

\begin{prop} Let $(p,q)\in\lc\times\lc$, we have
\begin{equation}\label{akns-est_kappa_petit o}
\fl\qquad\kappa_{a,n}(p,q)=\frac{(-1)^n}{\left(\left[|n|+\frac{a}{2}\right]\pi\right)^a}(1+o(1))=\frac{(-1)^n}{|n\pi|^a}(1+o(1))\quad,\quad
|n|\rightarrow+\infty.
\end{equation}
\end{prop}

\noindent\textit{Proof.} Introducing
\eqref{akns-estimation_L2_loc} in the $\kappa_{a,n}$ definition
leads to
$$\kappa_{a,n}=\frac{1}{{\lambda_{a,n}}^a}\Big(\ja{a-1}{_{a,n}}\cos{\beta}+\ja{a}{_{a,n}}\sin{\beta}+o(1)\Big).$$
Relation \eqref{annexe-bessel-ja-sinus} implies
$$\kappa_{a,n}=\frac{1}{{\lambda_{a,n}}^a}\Big(\cos{\left(\lambda_{a,n}-\frac{a\pi}{2}-\beta\right)}+o(1)\Big).$$
Now, with \eqref{akns-est_lambda_petit o}, we get
\begin{eqnarray*}
    \kappa_{a,n}&=&\frac{1}{{(n+\sgn{n}\frac{a}{2})}^a\pi^a}\left(\cos{\left(n\pi+(\sgn{(n)}-1)\frac{a\pi}{2}\right)}+o(1)\right),\\
    &=&\frac{(-1)^n}{{(n+\sgn{n}\frac{a}{2})}^a\pi^a}\left(\cos{\left[a\pi\frac{\sgn{(n)}-1}{2}\right]}+o(1)\right).
\end{eqnarray*}
Setting the signum of $n$ gives the result. \qed

\subsection{Orthogonality relations}\label{akns-section-ortho}

The following results, especially the corollary, confirm the
choice of the additional data: we have added only complementary
data. As in \cite{gg}, we obtain

\begin{prop} For all $(j,k)\in\Z^2$, we have
    \begin{enumerate}
        \item $\left\langle\nabla_{p,q}\lambda_{a,j},{\nabla_{p,q}\lambda_{a,k}}^\perp\right\rangle=0,$
        \item $\left\langle A_j(\cdot,p,q),{\nabla_{p,q}\lambda_{a,k}}^\perp\right\rangle=\delta_{j,k},$
        \item $\left\langle A_j(\cdot,p,q),{A_k(\cdot,p,q)}^\perp\right\rangle=0$.
    \end{enumerate}
\end{prop}
%\noindent\textit{Proof.} See \cite{moi-these} for details. \qed

Before giving the corollary, be more specific:
\begin{defn} A vector family $(u_k)_{k\in\Z}$ of an Hilbert space
is called free or its elements are linearly independent if each
element of the family is not in the closed span of the others.
More precisely:
$$\forall k\in\Z\,,\quad u_k\notin \overline{\vect{u_j|j\in\Z, j\neq k}}.$$
\end{defn}

\begin{coro}\label{akns-gradients_libres} For all $(j,k)\in\Z^2$, we have
\begin{enumerate}
    \item $\left\langle\nabla_{p,q}\kappa_{a,j},{\nabla_{p,q}\kappa_{a,k}}^\perp\right\rangle=0,$
    \item $\left\langle\nabla_{p,q}\kappa_{a,j},{\nabla_{p,q}\lambda_{a,k}}^\perp\right\rangle=\kappa_{a,j}(p,q)\delta_{j,k}.$
\end{enumerate}
$\left(\nabla_{p,q}\lambda_{a,n}\right)_{n\in\Z}\cup\left(\nabla_{p,q}\kappa_{a,n}\right)_{n\in\Z}$
is a free family in $\lr\times\lr$.
\end{coro}
%\noindent\textit{Proof.}
%    The first relations follow from the proposition and
%    equation \eqref{akns-expr_grad_kappa}. For the linear independence of the family,
%    consider for instance $\nabla_{p,q}\lambda_{a,k}$. Orthogonality
%    relations give on the one hand
%    $$\left\langle\nabla_{p,q}\lambda_{a,k},{\nabla_{p,q}\kappa_{a,k}}^\perp\right\rangle\neq
%    0$$
%    and on the other hand
%    \begin{eqnarray*}
%    &\left\langle\nabla_{p,q}\lambda_{a,j},{\nabla_{p,q}\kappa_{a,k}}^\perp\right\rangle=0,& \quad \mathrm{ for all }\,j\neq k,\\
%    &\left\langle\nabla_{p,q}\kappa_{a,j},{\nabla_{p,q}\kappa_{a,k}}^\perp\right\rangle=0,& \quad \mathrm{
%    for all  }\,j\in\Z.
%    \end{eqnarray*}
%    In other words
%    $$\left\{\begin{array}{c}
%      \nabla_{p,q}\lambda_{a,k}\notin\left(\vect{\nabla_{p,q}\kappa_{a,k}}\right)^\perp, \\
%      \\
%      \overline{\vect{\nabla_{p,q}\lambda_{a,j}, {j\neq k}; \nabla_{p,q}\kappa_{a,j},j\in\Z}}\subset
%      \left(\vect{\nabla_{p,q}\kappa_{a,k}}\right)^\perp. \\
%    \end{array}\right.$$
%     The same property can be obtained for $\nabla_{p,q}\kappa_{a,k}$, and thus we get the result.
%\qed

%   Construction de l'application spectrale

\subsection{The spectral map}

Introduce the quantities $\widetilde{\lambda}_ {a,n}(p,q)$ and
$\widetilde{\kappa}_{a,n}(p,q)$ such that
\begin{eqnarray*}
    \lambda_{a,n}(p,q)=\left(n+\sgn(n)\frac{a}{2}\right)\pi+\beta+\widetilde{\lambda}_{a,n}(p,q).\\
    \kappa_{a,n}(p,q)=\frac{(-1)^{n}}{\left[\left(|n|+\frac{a}{2}\right)\pi\right]^a}\left(1+\widetilde{\kappa}_{a,n}(p,q)\right).
\end{eqnarray*}
Now, with the estimates \eqref{akns-est_lambda_petit o} and
\eqref{akns-est_kappa_petit o}, we define the spectral map
$\lambda^a\times\kappa^a: \lr\times\lr\rightarrow c_0(\Z)\times
c_0(\Z)$ by
\begin{equation}\label{akns-def_lambda_x_kappa}
  \left[\lambda^a\times\kappa^a\right](p,q)=\Big((\widetilde{\lambda}_{a,n}(p,q))_{n\in\Z},
   (\widetilde{\kappa}_{a,n}(p,q))_{n\in\Z}\Big),
\end{equation}
where $c_0(\Z)$ is the space of sequences $(u_n)_{n\in\Z}$ which
tend to $0$ when $|n|\rightarrow \infty$.

Following \cite{ist} or \cite{gr}, to obtain regularity of
$\lambda^a\times\kappa^a$ from its components, some uniformity is
needed. To this end, we introduce some transformation operators.

     %%%%%%%%%%%%%%%%%%%%%%%%%%%%%%%%%%%%%%%%%%%%%%%%%%%%%%%%%%%%%%%%%%%%%%%%%
%%%%%%%%%%%%%%%%%%%%%%%%%%%%%%%%%%%%%%%%%%%%%%%%%%%%%%%%%%%%%%%%%%%%%%%%%
\subsection{Transformations operators}
%%%%%%%%%%%%%%%%%%%%%%%%%%%%%%%%%%%%%%%%%%%%%%%%%%%%%%%%%%%%%%%%%%%%%%%%%
%%%%%%%%%%%%%%%%%%%%%%%%%%%%%%%%%%%%%%%%%%%%%%%%%%%%%%%%%%%%%%%%%%%%%%%%%
Such operators were first introduced by Guillot and Ralston in
\cite{gr} for the inverse spectral problem of the radial
Schrödinger operator when $a=1$; then used and extended to any
integer $a$ by Rundell and Sacks in \cite{rusa} and by the present
author in \cite{moi-SL}.

We construct similar operator adapted to the AKNS operator. An
important difference, excepted the matrix form, is a better
structure of the converse operators compared to the Schrödinger
operator. These operators turn to be adapted to the spectral data,
since both vectors family corresponding to $\lambda^a$ and
$\kappa^a$ are well transformed. The proofs of the following
lemmas are similar to those in \cite{rusa}. The main tool is the
use of Bessel function's properties (for a detailed proof see
\cite{moi-these}). Now, give some notations.
\begin{nota} For all $n\in\N$, let $U_n$ and
$V_n$ be defined by
$$U_n(x)=\Bigg[\begin{array}{c}
  0 \\
  x^n \\
\end{array}\Bigg]\quad\mathrm{ and }\quad
V_n(x)=\Bigg[\begin{array}{c}
  x^n \\
  0 \\
\end{array}\Bigg]\quad x\in[0,1].$$
\end{nota}

\begin{lem}\label{akns_lemme_reduction_indice} For all $a\in\N$, let
\begin{eqnarray*}
    \begin{array}{cccc}
  S_{a+1}:&\lc\times\lc&\longrightarrow&\lc\times\lc\\
            &(p,q)                       &\longmapsto&\Big(S_{a,1}[p]\, ,S_{a,2}[q]
            \Big)\\
\end{array}\\
\begin{array}{lccl}
  \mathrm{with}&S_{a,1}[p](x)&=&\displaystyle p(x)-2(2a+1)x^{2a}\int_x^1
\frac{p(t)}{t^{2a+1}}dt,\\
  \mathrm{and}  &S_{a,2}[q](x)&=&\displaystyle q(x)-2(2a+1)x^{2a+1}\int_x^1
  \frac{q(t)}{t^{2a+2}}dt.
\end{array}
\end{eqnarray*}

Moreover, we set $S_0:=\id_{\lc\times\lc}$. We have the following
properties:
\begin{enumerate}[(i)]
%%%%%%%%%%%%%%%%%%%%%%%%%%%%%%%%%%%%%%%%%%%%%%%%%%%%%%%%%%%%%%%%%%%%%%%%%
    \item The adjoint of $S_{a+1}$ is $\displaystyle {S_{a+1}}^\ast[f,g]=\Big({S_{a,1}}^\ast[f]\, ,{S_{a,2}}^\ast[g]\Big)$
    where
    \begin{eqnarray*}
      {S_{a,1}}^\ast[f](x) &=& f(x)-\frac{2(2a+1)}{x^{2a+1}}\int_0^x
      t^{2a}f(t)dt, \\
      {S_{a,2}}^\ast[g](x) &=& g(x)-\frac{2(2a+1)}{x^{2a+2}}\int_0^x
      t^{2a+1}g(t)dt.
    \end{eqnarray*}
%%%%%%%%%%%%%%%%%%%%%%%%%%%%%%%%%%%%%%%%%%%%%%%%%%%%%%%%%%%%%%%%%%%%%%%%%
    \item \label{akns-commute} The family $\{S_a\}$ pairwise commutes:
    $S_{a}S_{b}=S_{b}S_{a}$ for all $(a,b)\in\N^2$.
%%%%%%%%%%%%%%%%%%%%%%%%%%%%%%%%%%%%%%%%%%%%%%%%%%%%%%%%%%%%%%%%%%%%%%%%%
    \item $S_a$ is bounded on $\lc\times\lc$.
%%%%%%%%%%%%%%%%%%%%%%%%%%%%%%%%%%%%%%%%%%%%%%%%%%%%%%%%%%%%%%%%%%%%%%%%%
    \item Let $N_{a+1}:=\ker {S_{a+1}}^\ast$, then
    $N_{a+1}=\mathrm{Vect}(U_{2a},V_{2a+1})$.
%%%%%%%%%%%%%%%%%%%%%%%%%%%%%%%%%%%%%%%%%%%%%%%%%%%%%%%%%%%%%%%%%%%%%%%%%
    \item $S_{a+1}$ is a linear isomorphism between $\lc\times\lc$
    and ${N_{a+1}}^\perp$.\\
    Its inverse is the bounded operator on $\lc\times\lc$ defined
    by
    $$A_{a+1}[f,g]:=\Big({S_{a,2}}^{\ast}[f]\, ,{S_{a,1}}^{\ast}[g]\Big).$$
%%%%%%%%%%%%%%%%%%%%%%%%%%%%%%%%%%%%%%%%%%%%%%%%%%%%%%%%%%%%%%%%%%%%%%%%%
    \item \label{akns-reduction}$\Phi_a$ and $\Psi_a$ defined by
    $$\Phi_a(x)=\Bigg[\begin{array}{c}
      -2 j_{a-1}(x)j_{a}(x) \\
      j_{a}(x)^2-j_{a-1}(x)^2 \\
    \end{array}\Bigg]$$
    and $$\Psi_a(x)=\Bigg[\begin{array}{c}
      -\eta_{a-1}(x)j_{a}(x)-\eta_a(x)j_{a-1}(x) \\
      -\eta_{a-1}(x)j_{a-1}(x)+\eta_{a}(x)j_{a}(x) \\
    \end{array}\Bigg]$$ satisfy the relations
    $$\Phi_{a+1}=-{S_{a+1}}^\ast[\Phi_{a}]\quad \mathrm{ and }\quad
    \Psi_{a+1}=-{S_{a+1}}^\ast[\Psi_{a}].$$
%%%%%%%%%%%%%%%%%%%%%%%%%%%%%%%%%%%%%%%%%%%%%%%%%%%%%%%%%%%%%%%%%%%%%%%%%
\end{enumerate}
\end{lem}

\begin{lem}\label{akns-akns_lemme_bessel_sinus}
For all $a\in\N$ we define $T_a$ by
\begin{equation}\label{akns-akns_eq_lemme_bessel_sinus_def}
T_a=(-1)^{a+1}S_a S_{a-1}\cdots S_1\quad,\quad T_0=-S_0.
\end{equation}
Let $T_a[f,g]=\Big({T_a^1}[f]\, ,{T_a^2}[g]\Big)$, then
\begin{enumerate}[(i)]
    \item $T_a$ is a bounded, one-to-one operator on $\lc\times\lc$
        such that for all $p,q\in\lc$ and all $\lambda\in\C$
        \begin{eqnarray}\label{akns-akns_eq_lemme_bessel_sinus_1}
            \fl\qquad \int_0^1 \Phi_a(\lambda t)\cdot\Bigg[\begin{array}{cc}
            p(t) \\  q(t) \\
        \end{array}\Bigg]dt=\int_0^1 \Bigg[\begin{array}{cc}
            \sin(2\lambda t) \\ \cos(2\lambda t) \\
        \end{array}\Bigg]\cdot T_a[p,q](t)dt,\\
        \label{akns-akns_eq_lemme_bessel_sinus_2}
            \fl\qquad\int_0^1 \Psi_a(\lambda t)\cdot\Bigg[\begin{array}{c}
            p(t) \\
            q(t) \\
        \end{array}\Bigg]dt=\int_0^1 \Bigg[\begin{array}{c}
        \cos(2\lambda t) \\
        -\sin(2\lambda t) \\
        \end{array}\Bigg]\cdot T_a[p,q](t)dt.
        \end{eqnarray}
    \item The adjoint of $T_a$, $\displaystyle T_{a}^\ast[f,g]=\Big({T_a^1}^\ast[f]\, ,{T_a^2}^\ast[g]\Big)$
    verifies
    \begin{equation}\label{akns-akns_eq_lemme_bessel_sinus_3}
    \fl\qquad \Phi_a(\lambda x)=T_a^\ast\Bigg[\begin{array}{c}
        \sin(2\lambda x) \\ \cos(2\lambda x) \\
    \end{array}\Bigg]\quad\mathrm{ and }\quad \Psi_a(\lambda x)=T_a^\ast\Bigg[\begin{array}{c}
        \cos(2\lambda x) \\ -\sin(2\lambda x) \\
    \end{array}\Bigg]
    \end{equation}
    and $$\mathrm{Ker}(T_a^\ast)=\displaystyle\bigoplus_{k=1}^a N_{k}.$$
    \item $T_a$ defines a linear isomorphism between $\lc\times\lc$
    and $\left(\displaystyle\bigoplus_{k=1}^a N_{k}\right)^\perp$.\\
    Its inverse is the bounded operator on $\lc\times\lc$ defined by
    $$B_a[f,g]:=\Big({T_a^2}^{\ast}[f]\, ,{T_a^1}^{\ast}[g]\Big).$$
\end{enumerate}
\end{lem}
%%%%%%%%%%%%%%%%%%%%%%%%%%%%%%%%%%%%%%%%%%%%%%%%%%%%%%%%%%%%%%%%%%%%%%%%%

     %   Asymptotiques

\subsection{Asymptotics upgrade}

%Nous allons maintenant combler le retard par rapport à l'opérateur
%de Schrödinger. En effet, nous sommes maintenant en mesure de
%préciser les asymptotiques des diverses grandeurs qui nous
%intéressent.

The following asymptotics are delicate to obtain since we want
them to figure both asymptotic behavior with respect to $n$ and
singular behavior with respect to $x$. Transformation operator
will help us to handle this difficulty.

First, give a tool ensuring us some uniformity with respect to
potentials. It is a Riemann-Lebesgue type lemma:

%%%%%%%%%%%%%%%%%%%%%%%%%%%%%%%%%%%%%%%%%%%%%%%%%%%%%%%%%%%%%%%%%%%%%%%%%%
\begin{lem}[Lemma A.1. in \cite{ag}(See also \cite{mis})]
\label{akns-lemme_AG_A1} $$\left(\int_0^1
f(t)e^{2\imath\pi(k+\varepsilon_k)t}dt\right)_{k\in\Z}\in\ell_\C^2(\Z)$$
uniformly with respect to
$\left(f,(\varepsilon_k)_{k\in\Z}\right)$ on bounded sets of
$\lc\times\ell_\C^\infty(\Z).$
\end{lem}
%%%%%%%%%%%%%%%%%%%%%%%%%%%%%%%%%%%%%%%%%%%%%%%%%%%%%%%%%%%%%%%%%%%%%%%%%%
Give an useful writing shortcut:
\begin{nota}
Let $(f_n)_{n\in\Z}$ a sequence of $L^\infty_\C(0,1)$ functions.
The equality
$$f_n(x)=\ell^2(n),\quad x\in[0,1], \quad n\in\Z$$
means
$$(\n{f_n}{\infty})_{n\in\Z}\in\ell^2_\R(\Z).$$
\end{nota}
%\begin{rem} Il est aussi possible de définir cette notation pour
%des suites de fonctions de $\lc$, mais celles qui apparaissent
%dans notre étude sont bornées et de plus la comparaison des normes
%suivantes
%$$\n{f}{\lc} \leq \n{f}{\infty},\quad f\in L^\infty_\C(0,1)$$
%permet de s'affranchir de cette précision.
%\end{rem}

%\subsubsection{Estimation de la solution régulière}

\begin{thm}\label{akns-theoreme_asymptotique_L2}
Uniformly on $[0,1]$ and locally uniformly on $\lc\times\lc$ we
have the following estimate:
\begin{equation}\label{akns-eq_asymptotique_sol_reg_L2}
\fl\qquad
\Big|\sr(x,\lambda_{a,n}(p,q),p,q)-R(x,\lambda_{a,n}(p,q))\Big|\leq
C \left[\frac{x}{1+|\lambda_{a,n}|
    x}\right]^a\ell^2(n),\quad |n|\rightarrow \infty,
\end{equation}
and locally uniformly on $\lc\times\lc$, we have
\begin{equation}\label{akns-eq_asymptotique_vp_L2}
    \lambda_{a,n}(p,q)=\left(n+\sgn(n)\frac{a}{2}\right)\pi+\beta+\ell^2(n),\quad |n|\rightarrow \infty.
\end{equation}
\end{thm}
%%%%%%%%%%%%%%%%%%%%%%%%%%%%%%%%%%%%%%%%%%%%%%%%%%%%%%%%%%%%%%%%%%%%%%%%%%
%%%%%%%%%%%%%%%%%%%%%%%%%%%%%%%%%%%%%%%%%%%%%%%%%%%%%%%%%%%%%%%%%%%%%%%%%%
\noindent\textit{Proof.} We first prove a similar estimate for
$R_1(x,\lambda_{a,n}(p,q),p,q)$. For this, recall (see the proof
of Theorem \ref{akns-thm_estimation_h1_reg}) that
$R_1(x,\lambda,p,q)=\lambda^{-a}[X(q)+Y(p)]$. Thus notations from
Lemma \ref{akns_lemme_reduction_indice} give
\begin{eqnarray*}
\fl\qquad
R_1(x,\lambda,p,q)=&\frac{1}{\lambda^a}\Bigg[\begin{array}{c}
  \displaystyle-\na{a-1}{x} \\
 \displaystyle \na{a}{x}\\
\end{array}\Bigg]\int_0^1 \Phi_a(\lambda t)\cdot
\Bigg[\begin{array}{c}
  \indic_{[0,x]}(t)p(t) \\
  \indic_{[0,x]}(t)q(t) \\
\end{array}\Bigg]dt\\&+\frac{1}{\lambda^a}\Bigg[\begin{array}{c}
\displaystyle\ja{a-1}{x}\\
\displaystyle-\ja{a}{x}\\
\end{array}\Bigg]\int_0^1 \Psi_a(\lambda t)\cdot
\Bigg[\begin{array}{c}
  \indic_{[0,x]}(t)p(t) \\
  \indic_{[0,x]}(t)q(t) \\
\end{array}\Bigg]dt.
\end{eqnarray*}
Estimate (\ref{akns-asymptotique_vp_L2_1}) implies that
$\lambda_{a,n}=n\pi+\varepsilon_n$ with
$(\varepsilon_n)_n\in\ell_\C^\infty(\Z)$. Then, lemmas
\ref{akns-akns_lemme_bessel_sinus} and \ref{akns-lemme_AG_A1} give
uniformly on $[0,1]$ and locally uniformly on $\lc\times\lc$:
\begin{eqnarray*}
  \left(\int_0^1 \Phi_a(\lambda_{a,n} t)\cdot \Bigg[\begin{array}{c}
  \indic_{[0,x]}(t)p(t) \\
  \indic_{[0,x]}(t)q(t) \\
\end{array}\Bigg]dt\right)_{n\in\Z}\in\ell_\C^2(\Z),\\
\left(\int_0^1 \Psi_a(\lambda_{a,n} t)\cdot \Bigg[\begin{array}{c}
  \indic_{[0,x]}(t)p(t) \\
  \indic_{[0,x]}(t)q(t) \\
\end{array}\Bigg]dt\right)_{n\in\Z}\in\ell_\C^2(\Z),
\end{eqnarray*}
in other words \begin{eqnarray*}
     \fl\qquad   R_1(x,\lambda_{a,n}(p,q),p,q)=\frac{\ell^2(n)}{{\lambda_{a,n}}^a}\Bigg[\begin{array}{c}
  \displaystyle-\na{a-1}{_{a,n} x} \\
 \displaystyle \na{a}{_{a,n} x}\\
\end{array}\Bigg]
+\frac{\ell^2(n)}{{\lambda_{a,n}}^a}\Bigg[\begin{array}{c}
\displaystyle\ja{a-1}{_{a,n} x}\\
\displaystyle-\ja{a}{_{a,n} x}\\
\end{array}\Bigg].
\end{eqnarray*} From \eqref{annexe-bessel_est_ja}, we obtain
\begin{equation}\label{akns-preuve_th_asym_eq1}
\left|\frac{\ell^2(n)}{{\lambda_{a,n}}^a}\Bigg[\begin{array}{c}
\displaystyle\ja{a-1}{_{a,n} x}\\
\displaystyle-\ja{a}{_{a,n} x}\\
\end{array}\Bigg]\right|\leq
C\left(\frac{x}{1+|\lambda_{a,n}|x}\right)^a\ell^2(n).
\end{equation}
For the first term in $R_1(x,\lambda_{a,n}(p,q),p,q)$ we split
$[0,1]$ in two:
\begin{description}
    \item[$|\lambda_{a,n} x|\geq 1$:] Since uniformly on $[0,1]$,
            $$\displaystyle\int_0^1 \Phi_a(\lambda_{a,n} t)\cdot
        \Bigg[\begin{array}{c}
        \indic_{[0,x]}(t)p(t) \\
        \indic_{[0,x]}(t)q(t) \\
        \end{array}\Bigg]dt=\ell^2(n)$$ and
        $$1=\frac{1+|\lambda_{a,n}|x}{1+|\lambda_{a,n}|x}\leq \frac{2|\lambda_{a,n}|x}{1+|\lambda_{a,n}|x}
        ,$$ we get
        $$\left|\int_0^1 \Phi_a(\lambda_{a,n} t)\cdot
        \Bigg[\begin{array}{c}
        \indic_{[0,x]}(t)p(t) \\
        \indic_{[0,x]}(t)q(t) \\
        \end{array}\Bigg]dt\right|\leq  \left(\frac{2|\lambda_{a,n}|x}{1+|\lambda_{a,n}|x}\right)^{2a}\ell^2(n).$$
    \item[$|\lambda_{a,n}x|\leq 1$:] Estimate
    \eqref{annexe-bessel_est_ja} gives
        $$\left|\int_0^1 \Phi_a(\lambda_{a,n} t)\cdot
        \Bigg[\begin{array}{c}
        \indic_{[0,x]}(t)p(t) \\
        \indic_{[0,x]}(t)q(t) \\
        \end{array}\Bigg]dt\right|\leq  C \left(\frac{|\lambda_{a,n}|x}{1+|\lambda_{a,n}|x}\right)^{2a}\int_0^x \Bigg[\begin{array}{c}
        |p(t)| \\
        |q(t)| \\
        \end{array}\Bigg]dt,$$
        where $C>0$ is uniform in $x$ and $n$, then
        $$\left|\int_0^1 \Phi_a(\lambda_{a,n} t)\cdot
        \Bigg[\begin{array}{c}
        \indic_{[0,x]}(t)p(t) \\
        \indic_{[0,x]}(t)q(t) \\
        \end{array}\Bigg]dt\right|\leq  C \negthinspace \left(\frac{|\lambda_{a,n}|x}{1+|\lambda_{a,n}|x}\right)^{2a}\negthinspace\negthinspace\int_0^{|\lambda_{a,n}|^{-1}}\negthinspace \Bigg[\begin{array}{c}
        |p(t)| \\
        |q(t)| \\
        \end{array}\Bigg]dt.$$
        Lemma \ref{annexe-lemme-carlson-1} gives the good bound.
\end{description}
Combining theses two estimates, we get uniformly on $[0,1]$ and
locally uniformly on $\lc\times\lc$:
\begin{equation}\label{akns-preuve_th_asym_eq2}
    \left|\int_0^1 \Phi_a(\lambda_{a,n} t)\cdot
        \Bigg[\begin{array}{c}
        \indic_{[0,x]}(t)p(t) \\
        \indic_{[0,x]}(t)q(t) \\
        \end{array}\Bigg]dt\right|\leq  C' \left(\frac{|\lambda_{a,n}|x}{1+|\lambda_{a,n}|x}\right)^{2a}\ell^2(n).
\end{equation}
Estimate \eqref{annexe-bessel_est_na} together with
\eqref{akns-preuve_th_asym_eq1} and
\eqref{akns-preuve_th_asym_eq2} gives
$$\Big|R_1(x,\lambda_{a,n}(p,q),p,q)\Big|\leq \left(\frac{x}{1+|\lambda_{a,n}|x}\right)^{a}\ell^2(n)$$
locally uniformly on $\lc\times\lc$ and uniformly on $[0,1]$.

With the recurrence relation and the estimation for
$\G(x,t,\lambda)$, follows uniformly on $[0,1]$ and locally
uniformly on $\lc\times\lc$:
\begin{eqnarray*}
\fl\qquad\Big|R_{k+1}(x,\lambda_{a,n}(p,q),p,q)\Big|\leq
\frac{C^k}{k!} \left(\int_0^x\left(|p(t)|+|q(t)|\right)dt\right)^k
\left(\frac{x}{1+|\lambda_{a,n}|x}\right)^{a}\ell^2(n),
\end{eqnarray*} summing up, we get the result.
Eigenvalues estimate is deduced directly from
$\sr(x,\lambda_{a,n}(p,q),p,q)$'s estimate and from
\eqref{akns-asymptotique_vp_L2_1}-\eqref{akns-asymptotique_vp_L2_2}.
\qed

     %%%%%%%%%%%%%%%%%%%%%%%%%%%%%%%%%%%%%%%%%%%%%%%%%%%%%%%%%%%%%%%%%%%%%%%%%%
%%%%%%%%%%%%%%%%%%%%%%%%%%%%%%%%%%%%%%%%%%%%%%%%%%%%%%%%%%%%%%%%%%%%%%%%%%

In a very similar way, we upgrade the control of the singular
solution and doing it justify the choice and existence of the
singular solution as announced in the first remark.

%%%%%%%%%%%%%%%%%%%%%%%%%%%%%%%%%%%%%%%%%%%%%%%%%%%%%%%%%%%%%%%%%%%%%%%%%%
%%%%%%%%%%%%%%%%%%%%%%%%%%%%%%%%%%%%%%%%%%%%%%%%%%%%%%%%%%%%%%%%%%%%%%%%%%
\begin{thm}\label{akns-theoreme_asymptotique_L2_sing}
Let $(p,q)\in\lc\times\lc$, then uniformly on $(0,1]$ and locally
uniformly on $\lc\times\lc$ we have:
\begin{equation}\label{akns-eq_asymptotique_sol_sing_L2}
    \Big|\ss(x,\lambda_{a,n}(p,q),p,q)-S(x,\lambda_{a,n}(p,q))\Big|\leq C\left[\frac{1+|\lambda_{a,n}|
    x}{x}\right]^a\ell^2(n).
\end{equation}
\end{thm}
%%%%%%%%%%%%%%%%%%%%%%%%%%%%%%%%%%%%%%%%%%%%%%%%%%%%%%%%%%%%%%%%%%%%%%%%%%
\noindent\textit{Proof.}
As for the regular solution, we obtain (see \cite{moi-these}) the
uniform estimate in $x\in[0,1]$ and locally uniform on
$\lr\times\lr$:
$$\tilde{\ss}(x,\lambda_{a,n}(p,q),p,q)=S(x,\lambda)+\mathcal{O}\left(\left[\frac{1+|\lambda_{a,n}|
    x}{x}\right]^a\right)\ell^2(n).$$
Then, we get easily $\wr(\lambda_{a,n}(p,q),p,q) =
\wr(\lambda_{a,n}(p,q),0)+\ell^2(n)= 1+\ell^2(n)$ and through
$\ss(x,\lambda,p,q)=\displaystyle\frac{\tilde{\ss}(x,\lambda,p,q)}{\wr(\lambda,p,q)}$,
we reach the result. \qed

%%%%%%%%%%%%%%%%%%%%%%%%%%%%%%%%%%%%%%%%%%%%%%%%%%%%%%%%%%%%%%%%%%%%%%

Straightforward calculations let us deduce the following
estimations:
%%%%%%%%%%%%%%%%%%%%%%%%%%%%%%%%%%%%%%%%%%%%%%%%%%%%%%%%%%%%%%%%%%%%
\begin{coro}\label{akns-thm_asympt_misc} Uniformly on $[0,1]$ and locally
uniformly on $\lr\times\lr$, when $|n|\rightarrow\infty$, we have
\begin{eqnarray}
% Ajout pour l'article en anglais
    \label{akns-est-norm-Rn} \|\sr_n(\cdot,p,q)\|^2 = \frac{1}{{\lambda_{a,n}}^{2a}}\left(1+
    \ell^2(n)\right),\\
    \label{akns-produit-scalaire}
                \Big\langle \sr_n(\cdot,p,q),\ss_n(\cdot,p,q)\Big\rangle
                =\ell^2(n),\\
% Fin de l'ajout
    \label{akns-eq_asympt_misc_1}G_n(x,p,q)=\Bigg[\begin{array}{c}\ja{a-1}{_{a,n} x} \\
        -\ja{a}{_{a,n} x} \\\end{array}\Bigg]+\ell^2(n),\\
    \label{akns-eq_asympt_misc_2}\nabla_{p,q}\lambda_{a,n}(p,q)=\Phi_a(\lambda_{a,n} x)+\ell^2(n),\\
    \label{akns-eq_asympt_misc_3}\kappa_{a,n}(p,q)=\frac{(-1)^{n}}{\left[\left(|n|+\frac{a}{2}\right)\pi\right]^a}\left[1+\ell^2(n)\right]=
    \frac{(-1)^{n}}{|n\pi|^a}\left[1+\ell^2(n)\right],\\
    \label{akns-eq_asympt_misc_4} A_n(x,p,q)=\Psi_a(\lambda_{a,n} x)+\ell^2(n),\\
    \label{akns-eq_asympt_misc_5}\frac{\nabla_{p,q}\kappa_{a,n}(p,q)}{\kappa_{a,n}(p,q)}=\Psi_a(\lambda_{a,n} x)+\ell^2(n).
\end{eqnarray}
\end{coro}
%%%%%%%%%%%%%%%%%%%%%%%%%%%%%%%%%%%%%%%%%%%%%%%%%%%%%%%%%%%%%%%%%%%%
Now, the spectral map can be correctly defined by
$$
\begin{array}{cccc}
  \lambda^a\times\kappa^a: & \lr\times\lr & \longrightarrow & \ell^2(\Z)\times\ell^2(\Z)\\
   & (p,q) & \longmapsto & \Big((\widetilde{\lambda}_{a,n}(p,q))_{n\in\Z},
   (\widetilde{\kappa}_{a,n}(p,q))_{n\in\Z}\Big), \\
\end{array}
$$
and, following \cite{ist} and  \cite{gr}, previous analyticity
results and the local uniformity with respect to the potentials
give us:

\begin{thm}\label{akns-appl-analyt} $\lambda^a\times\kappa^a$ is a real-analytic
map on $\lr\times\lr$.\\
     Its Fréchet derivative is given by the linear map from
    $\lr\times\lr$ to $\ell^2(\Z)\times\ell^2(\Z)$:
    $$d_{p,q}(\lambda^a\times\kappa^a)(v)=\Big(\left(\left\langle\nabla_{p,q}\lambda_{a,n},v\right\rangle\right)_{n\in\Z},
    \left(\left\langle\nabla_{p,q}\widetilde{\kappa}_{a,n},v\right\rangle\right)_{n\in\Z}\Big).$$
\end{thm}

     %%%%%%%%%%%%%%%%%%%%%%%%%%%%%%%%%%%%%%%%%%%%%%%%%%%%%%%%%%%%%%%%%%%%%%%%%%%
\section{The inverse spectral problem}
%%%%%%%%%%%%%%%%%%%%%%%%%%%%%%%%%%%%%%%%%%%%%%%%%%%%%%%%%%%%%%%%%%%

Now, give the main result

\begin{thm}\label{akns-akns_th_diff_inversible}$\quad$\\
    $d_{p,q}(\lambda^a\times\kappa^a)$ is an isomorphism between $\lr\times\lr$
    and
    $\ell^2(\Z)\times\ell^2(\Z)$.
\end{thm}
\noindent\textit{Proof.} In view of the relation
$$\nabla_{p,q}\widetilde{\kappa}_{a,n}=(-1)^n\left[\left(|n|+\frac{a}{2}\right)\pi\right]^a\nabla_{p,q}\kappa_{a,n},$$
corollary \ref{akns-gradients_libres} implies that
$\left(\nabla_{p,q}\lambda_{a,n}\right)_{n\in\Z}\cup
\left(\nabla_{p,q}\widetilde{\kappa}_{a,n}\right)_{n\in\Z}$ is a
free family in $\lr\times\lr$. Let define $r_n$ and $s_n$ by
\begin{eqnarray}
\label{akns-def-rn}  r_n(x)=\nabla_{p,q}\lambda_{a,n}(x)-\Phi_a(\lambda_{a,n} x),\\
\label{akns-def-sn}
s_n(x)=\nabla_{p,q}\widetilde{\kappa}_{a,n}(x)-\Psi_a(\lambda_{a,n}
x).
\end{eqnarray}
With lemma \ref{akns-akns_lemme_bessel_sinus}, we have for all
$v\in\lr\times\lr$,
\begin{eqnarray}
    \label{akns-akns_eq_grad_pert_1}\left\langle\nabla_{(p,q)}\lambda_{a,n}(V),v\right\rangle=\int_0^1
    \Bigg(\Bigg[\begin{array}{c}
        \sin{(2\lambda_{a,n} t)} \\
        \cos{(2\lambda_{a,n} t)}\\
    \end{array}\Bigg]+R_n(t)\Bigg)\cdot T_a[v](t)dt,\\
    \label{akns-akns_eq_grad_pert_2}\left\langle\nabla_{(p,q)}\widetilde{\kappa}_{a,n}(V),v\right\rangle=\int_0^1
    \Bigg(\Bigg[\begin{array}{c}
        \cos{(2\lambda_{a,n} t)} \\
        -\sin{(2\lambda_{a,n} t)}\\
    \end{array}\Bigg]+S_n(t)\Bigg)\cdot T_a[v](t)dt,
\end{eqnarray}
where $R_n=B_a^\ast[r_n]$ and $S_n=B_a^\ast[s_n]$. Introduce
operator $F$ defined by
\begin{eqnarray*}
\fl
%\qquad
        F(w)=\Bigg(\Bigg\{
        \Bigg\langle\Bigg[\begin{array}{c} \sin{(2\lambda_{a,n} t)} \\
        \cos{(2\lambda_{a,n} t)}\\\end{array}\Bigg]+R_n(t),w\Bigg\rangle
                \Bigg\}_{n\in\Z},%\\
                 %\qquad\qquad
                 \Bigg\{
        \Bigg\langle\Bigg[\begin{array}{c} \cos{(2\lambda_{a,n} t)} \\
        -\sin{(2\lambda_{a,n} t)}\\ \end{array}\Bigg]+S_n(t),w\Bigg\rangle
                \Bigg\}_{n\in\Z}\Bigg),
\end{eqnarray*}
in order to get $d_{p,q}(\lambda^a\times\kappa^a)(v)=F\circ
T_a[v]$. From lemma \ref{akns-akns_lemme_bessel_sinus}, $T_a$ is a
bijection between $\lc\times\lc$ and
$\left(\displaystyle\bigoplus_{k=1}^a N_{k}\right)^\perp$. Thus,
we have to prove that $F$ is a bijection between
$\left(\displaystyle\bigoplus_{k=1}^a N_{k}\right)^\perp$ and
$\ell^2(\Z)\times\ell^2(\Z)$. To this end, we will show that the
operator $\bf{F}$ sending functions in $\lr\times\lr$ into their
Fourier coefficients (or, in other words, the scalar products)
with respect to the family
\begin{eqnarray}
\fl\qquad        \mathcal{F}=\Bigg(\left\{
            U_{2k}\right\}_{k=0}^{a-1},&\Bigg\{
            \Bigg[\begin{array}{c} \sin{(2\lambda_{a,n} t)} \\
            \cos{(2\lambda_{a,n}
            t)}\\\end{array}\Bigg]+R_n(t)\Bigg\}_{n\in\Z},\nonumber\\
            &\left\{V_{2k+1}\right\}_{k=0}^{a-1},\Bigg\{
        \Bigg[\begin{array}{c} \cos{(2\lambda_{a,n} t)} \\
        -\sin{(2\lambda_{a,n} t)}\\
        \end{array}\Bigg]+S_n(t)\Bigg\}_{n\in\Z}\Bigg),\label{akns-def_famille_libre}
\end{eqnarray}
is a invertible map from $\lr\times\lr$ to
$\ell^2(\Z)\times\ell^2(\Z).$ For this, recall the following
property (see \cite{ist}: Appendix D, theorem 3).
%%%%%%%%%%%%%%%%%%%%%%%%%%%%%%%%%%%%%%%%%%%%%%%%%%%%%%%%%%%%%%%%%%%%%
\begin{lem}\label{lemme-pt} Let $\{f_n\}_{n\in\Z}$ be a free family of vectors in an Hilbert space
$H$ close to an orthonormal basis $\{e_n\}_{n\in\Z}$ of $H$, ie
$\sum \n{f_n-e_n}{2}^2<\infty$.\\
 Then
$\left\{f_n\right\}_{n\in\Z}$ is a basis for $H$ and the map
$\mathbf{F}:x\mapsto\{(f_n,x)\}_{n\in\Z}$ is a linear isomorphism
from $H$ onto $\ell^2(\Z)$.
\end{lem}
%%%%%%%%%%%%%%%%%%%%%%%%%%%%%%%%%%%%%%%%%%%%%%%%%%%%%%%%%%%%%%%%%%%%%
Estimates \eqref{akns-eq_asympt_misc_2},
\eqref{akns-eq_asympt_misc_3} and \eqref{akns-eq_asympt_misc_5}
lead to $r_n=\ell^2(n)$ and $s_n=\ell^2(n)$. Boundedness of
$B_a^\ast$ thus gives $R_n=\ell^2(n)$ and $S_n=\ell^2(n)$ which, together
with the orthogonal basis of $\lr\times\lr$ %(voir la remarque
%\ref{akns-base_amour} ci-dessous)
\begin{equation}\label{akns-exp-F0}
\fl\qquad    \mathcal{F}_0=\Bigg\{
        \Bigg[\begin{array}{c} \sin{\left(2\left((n+\frac{a}{2})\pi+\beta\right)t\right)} \\
        \cos{\left(2\left((n+\frac{a}{2})\pi+\beta\right)t\right)}\\\end{array}\Bigg],
        \Bigg[\begin{array}{c} \cos{\left(2\left((n+\frac{a}{2})\pi+\beta\right)t\right)} \\
        -\sin{\left(2\left((n+\frac{a}{2})\pi+\beta\right)t\right)}\\ \end{array}\Bigg],n\in\Z\Bigg\},
\end{equation}
and a correct arrangement of each vectors family (see remark
bellow %\ref{rem-numerotation}
), prove the closeness of
$\mathcal{F}$ and $\mathcal{F}_0$. Lemma
\ref{akns-akns_lemme_famille_libre} gives the freedom of
$\mathcal{F}$ and thus lemma \ref{lemme-pt} is applicable. \qed

%%%%%%%%%%%%%%%%%%%%%%%%%%%%%%%%%%%%%%%%%%%%%%%%%%%%%%%%%%%%%%%%%%%%

\noindent\textit{Remark.} At first sight, the ``loss" of
eigenvalues appeared in the counting lemma and the non-zero kernel
of the transformation operator seem to be barriers to solve the
inverse problem. In fact, it is not, it helps us to fit correctly
vectors family $\mathcal{F}$ and $\mathcal{F}_0$. Be more
specific: let $f_{n,1}^0$ and $f_{n,2}^0$ be defined by
\eqref{akns-exp-F0}, in other words, we just write
$\mathcal{F}_0=\left\{f_{n,1}^0,\, f_{n,2}^0,\, n\in\Z\right\}$.
For $\mathcal{F}$ we choose the following numbering: set
$\mathcal{F}=\left\{f_{n,1},\, f_{n,2},\, n\in\Z\right\}$ where
for any integer $n\geq 0$,
$$f_{n,1}(t)=\Bigg[\begin{array}{c} \sin{(2\lambda_{a,n} t)} \\
            \cos{(2\lambda_{a,n}
            t)}\\\end{array}\Bigg]+R_n(t),\quad f_{n,2}(t)=\Bigg[\begin{array}{c} \cos{(2\lambda_{a,n} t)} \\
        -\sin{(2\lambda_{a,n} t)}\\
        \end{array}\Bigg]+S_n(t),$$
for any integer $n$ such that $n\in\ing-a,-1\ind$,
$$f_{n,1}=U_{-2n-2},\quad f_{n,2}=V_{-2n-1},$$
and for all integer $n$ such that $n\leq -a-1$,
$$f_{n,1}(t)=\Bigg[\begin{array}{c} \sin{(2\lambda_{a,n+a} t)} \\
            \cos{(2\lambda_{a,n+a}
            t)}\\\end{array}\Bigg]+R_{n+a}(t), \, f_{n,2}(t)=\Bigg[\begin{array}{c} \cos{(2\lambda_{a,n+a} t)} \\
        -\sin{(2\lambda_{a,n+a} t)}\\
        \end{array}\Bigg]+S_{n+a}(t).$$
With this notation and using the eigenvalue estimate
\eqref{akns-eq_asymptotique_vp_L2}, for $j=1,2$, $(f_{n,j})_n$ is
asymptotically $\ell^2$-close to $(f_{n,j}^0)_n$ whenever
$n\rightarrow \pm\infty$.

In order to prove the freedom of $\mathcal{F}$, give a little
extension with the following
%%%%%%%%%%%%%%%%%%%%%%%%%%%%%%%%%%%%%%%%%%%%%%%%%%%%%%%%%%%%%%%%%%%%
\begin{prop}\label{prop-base}
Let $(E_{n,1},E_{n,2})_{n\in\Z}$ be a free vector family in
$\lr\times\lr$ satisfying the following properties:
\begin{enumerate}[(i)]
    \item Duality : there exists a bounded vector family $(F_{n,1},F_{n,2})_{n\in\Z}$ in $\lr\times\lr$,
    such that
    $$\langle E_{n,j},F_{m,j}\rangle=0,\quad (n,m)\in\Z^2,\quad j=1,2.$$
    $$\langle E_{n,1},F_{m,2}\rangle=\langle E_{n,2},F_{m,1}\rangle=\delta_{n,m},\quad \forall(n,m)\in\Z^2.$$
    \item Asymptotics:
    $$E_{n,1}=T_a^\ast\Bigg(\Bigg[\begin{array}{c} \sin{(2\lambda_{a,n} t)} \\
            \cos{(2\lambda_{a,n}
            t)}\\\end{array}\Bigg]+e_{n,1}\Bigg),\, E_{n,2}=T_a^\ast\Bigg(\Bigg[\begin{array}{c} \cos{(2\lambda_{a,n} t)} \\
     -\sin{(2\lambda_{a,n} t)}\\\end{array}\Bigg]+e_{n,2}\Bigg)$$
     with $\Big(\n{e_{n,j}}{\lr\times\lr}\Big)_{n\in\Z}\in\ell^2_\R(\Z)$, $j=1,2$.
    \item Summability: for any $k\in\ing
     0,2a-1\ind$, there exists $\omega\in\mathcal{C}_0^\infty([0,1],\R^2)$ such that
     for all $m\in\ing
     0,2a-1\ind $, $\left\langle\omega, W_m\right\rangle=\delta_{k,m}$ and
    $$\left(\langle \omega,e_{n,j}\rangle\right)_{n\in\Z}\in \ell^1_\R(\Z),\quad j=1,2.$$
\end{enumerate}
Then, the following family is free in $\lr\times\lr$
\begin{eqnarray*}
\fl\qquad        \mathcal{F}=\Bigg(\left\{
            U_{2k}\right\}_{k=0}^{a-1},&\Bigg\{
            \Bigg[\begin{array}{c} \sin{(2\lambda_{a,n} t)} \\
            \cos{(2\lambda_{a,n}
            t)}\\\end{array}\Bigg]+e_{n,1}(t)\Bigg\}_{n\in\Z},\\
            &\left\{V_{2k+1}\right\}_{k=0}^{a-1},\Bigg\{
        \Bigg[\begin{array}{c} \cos{(2\lambda_{a,n} t)} \\
        -\sin{(2\lambda_{a,n} t)}\\
        \end{array}\Bigg]+e_{n,2}(t)\Bigg\}_{n\in\Z}\Bigg)%,
        .
\end{eqnarray*}
%   VERSION ALTERNATIVE
%
%et
%\begin{equation*}
%        \mathcal{E}=\left\{E_{n,1}, E_{n,2}\right\}_{n\in\Z}.
%\end{equation*}
\end{prop}

\noindent\textit{Proof.} Since $T_a^\ast$ is bounded and
$\left(E_{n,1},E_{n,2}\right)_{n\in\Z}$ is free, condition
\textit{(ii)} implies the freedom of the following family
$$\Bigg\{\Bigg[\begin{array}{c} \sin{(2\lambda_{a,n} t)} \\
            \cos{(2\lambda_{a,n}
            t)}\\\end{array}\Bigg]+e_{n,1}(t)\Bigg\}_{n\in\Z}\cup\Bigg\{\Bigg[\begin{array}{c} \cos{(2\lambda_{a,n} t)} \\
     -\sin{(2\lambda_{a,n} t)}\\\end{array}\Bigg]+e_{n,1}(t)\Bigg\}_{n\in\Z}.$$

Let $k\in \ing 0,2a-1\ind$, we define $W_k$ by $W_{k}=U_{k}$ if
$k$ is even and $W_{k}=V_{k}$ otherwise. Show that $W_k$ is not in
the closure of
$\mathrm{Vect}\left(\mathcal{F}\setminus\left\{W_k\right\}\right)$.
(Precisely, we should prove iteratively that
$W_{k}\notin\overline{\vect{\mathcal{F}\setminus\left\{W_j,j\in\ing
k,2a-1\ind\right\}}}$, which is not necessary since it suffices to
set $\alpha_m^{(j)}=0$ for $m\in\ing k,2a-1\ind$ in the next
expression.) For this, suppose the contrary: there exists a vector
sequence defined for $j\in\N$ by
\begin{eqnarray*}
\fl\qquad    W_k^{(j)}(t)=&\sum_{m\in\ing 0,2a-1\ind, m\neq k}
\alpha_m^{(j)}\,W_m(t)+\sum_{n\in\ing-N_j,N_j\ind}a_n^{(j)}\Bigg(
            \Bigg[\begin{array}{c} \sin{(2\lambda_{a,n} t)} \\
            \cos{(2\lambda_{a,n}
            t)}\\\end{array}\Bigg]+e_{n,1}(t)\Bigg) \\
    &+\sum_{n\in\ing-N_j,N_j\ind}b_n^{(j)}\Bigg(\Bigg[\begin{array}{c} \cos{(2\lambda_{a,n} t)} \\
     -\sin{(2\lambda_{a,n} t)}\\\end{array}\Bigg]+e_{n,2}(t)\Bigg),
\end{eqnarray*}
with $N_j<\infty$, $ \alpha_m^{(j)},a_n^{(j)},b_n^{(j)}\in\R$ such
that $W_k^{(j)}\underset{j\rightarrow\infty}\longrightarrow W_k$
in $\lr\times\lr$. Recall that $T_a^\ast(W_m)=0$ for
$m=0,\dots,2a-1$, thus the sequence
$$w^{(j)}:=T_a^\ast(W_k^{(j)})=\sum_{n\in\ing-N_j,N_j\ind}a_n^{(j)}E_{n,1}+b_n^{(j)}E_{n,2}$$
converges towards $0$ in $\lr\times\lr$ when $j\rightarrow\infty$,
and point \textit{(i)} leads to
\begin{eqnarray}
    a_n^{(j)}=\int_0^1 w^{(j)}\cdot
    F_{n,2}\, d t\,\underset{j\rightarrow\infty}{\longrightarrow}0,\\
    b_n^{(j)}=\int_0^1 w^{(j)}\cdot
    F_{n,1}\, d t\,\underset{j\rightarrow\infty}{\longrightarrow}0.
\end{eqnarray}
and gives the uniform boundedness of $(a_n^{(j)})$ and
$(b_n^{(j)})$ with respect to $n$ and $j$.

Now consider %(voir l'Annexe \ref{annexe-ortho})
$\omega\in \classe{\infty}_0([0,1],\R^2)$ as in \textit{(iii)}.
Its smoothness and support property imply that for all $N\in\N$,
$$\int_0^1 \omega(t)\cdot \Bigg[\begin{array}{c} \sin{(2\lambda_{a,n} t)} \\
            \cos{(2\lambda_{a,n}
            t)}\\\end{array}\Bigg]d t,\int_0^1 \omega(t)\cdot \Bigg[\begin{array}{c} \cos{(2\lambda_{a,n} t)} \\
     -\sin{(2\lambda_{a,n} t)}\\\end{array}\Bigg]d t=\mathcal{O}\left(\frac{1}{n^N}\right).$$
Thus, second part of \textit{(iii)} shows the summability of
$$\Bigg\{\Bigg\langle\omega,t\mapsto \Bigg[\begin{array}{c} \sin{(2\lambda_{a,n} t)} \\
            \cos{(2\lambda_{a,n}
            t)}\\\end{array}\Bigg]+e_{n,1}(t)\Bigg\rangle\Bigg\}_{n\in\Z}
            $$
and
$$\Bigg\{\Bigg\langle\omega,t\mapsto \Bigg[\begin{array}{c} \cos{(2\lambda_{a,n} t)} \\
            -\sin{(2\lambda_{a,n}
            t)}\\\end{array}\Bigg]+e_{n,2}(t)\Bigg\rangle\Bigg\}_{n\in\Z}.$$
We complete the proof writing
\begin{eqnarray*} \fl\qquad\left\langle\omega,
W_k^{(j)}\right\rangle=&\sum_{n\in\ing-N_j,N_j\ind}a_n^{(j)}\Bigg\langle\omega,\Bigg(
    \Bigg[\begin{array}{c} \sin{(2\lambda_{a,n} t)} \\
    \cos{(2\lambda_{a,n}
    t)}\\\end{array}\Bigg]+e_{n,1}(t)\Bigg)\Bigg\rangle\\
    &+\sum_{n\in\ing-N_j,N_j\ind}b_n^{(j)}\Bigg\langle\omega,\Bigg(\Bigg[\begin{array}{c} \cos{(2\lambda_{a,n} t)} \\
     -\sin{(2\lambda_{a,n}
     t)}\\\end{array}\Bigg]+e_{n,2}(t)\Bigg)\Bigg\rangle,
\end{eqnarray*}
indeed, this shows, by dominated convergence, that %puis en appliquant le théorème de
%convergence dominée qui montre que
$$\left\langle\omega,
W_k^{(j)}\right\rangle\underset{j\rightarrow\infty}\longrightarrow
0,$$ which is in contradiction with the definition of $\omega$. So
$\mathcal{F}$ is a free family. \qed
%%%%%%%%%%%%%%%%%%%%%%%%%%%%%%%%%%%%%%%%%%%%%%%%%%%%%%%%%%%%%%%%%%%%

\begin{lem} \label{akns-akns_lemme_famille_libre}
$\mathcal{F}$ is a free family in $\lr\times\lr$.
\end{lem}

\noindent\textit{Proof.} Let us apply proposition \ref{prop-base}.
For this, we consider the following vectors
\begin{eqnarray*}
    E_{n,1}=\nabla_{p,q}\lambda_{a,n},\quad
    E_{n,2}=\nabla_{p,q}\widetilde{\kappa}_{a,n},\quad n\in\Z,\\
    F_{n,1}={\nabla_{p,q}\lambda_{a,n}}^\perp,\quad
    F_{n,2}=-{\nabla_{p,q}\widetilde{\kappa}_{a,n}}^\perp,\quad n\in\Z.
\end{eqnarray*}
Results from section \ref{akns-section-ortho} show that
$(E_{n,1},E_{n,2})_{n\in\Z}$ are linearly independent and that
condition \textit{(i)} is verified.

Relations \eqref{akns-akns_eq_grad_pert_1},
\eqref{akns-akns_eq_grad_pert_2} and estimates
\eqref{akns-eq_asympt_misc_2}, \eqref{akns-eq_asympt_misc_5} give
us condition \textit{(ii)} with
$$e_{n,1}={B_a}^\ast[r_n],\quad e_{n,2}={B_a}^\ast[s_n],$$
where $r_n$ and $s_n$ are defined by \eqref{akns-def-rn} and
\eqref{akns-def-sn}.

Now, condition \textit{(iii)} is left to be proved.\\
First, there exists $\omega\in \classe{\infty}_0([0,1],\R^2)$
compactly supported in $[\delta,1]$ for some $\delta>0$, such that
$m\in\ing 0,2a-1\ind $, $\left\langle\omega,
W_m\right\rangle=\delta_{k,m}$. Second, from the definition of
$S_a^\ast$ given in lemma \ref{akns_lemme_reduction_indice},
$B_a[\omega]$ is in $\classe{\infty}([0,1],\R^2)$ and supported in
$[\delta,1]$. We are now able to prove the summation properties.

Let $\varepsilon_n=(\varepsilon_n^1,\varepsilon_n^2)$ be defined
by $\varepsilon_n(x,V)=\sr_n(x,V)-R(x,\lambda_{a,n}(V))$ and plug
it in $\nabla_{p,q}\lambda_{a,n}$ via
\eqref{akns-expr_grad_lambda}. We get
\begin{eqnarray*}
\fl\quad  2G_{n,1}(x,V)G_{n,2}(x,V) &=
2\big(R_1(x,\lambda_{a,n})+\varepsilon_n^1\big)\big(R_2(x,\lambda_{a,n})+
    \varepsilon_n^2\big)\n{\sr_n(\cdot,p,q)}{2}^{-2}, \\
   &= \Big(2 R_1(x,\lambda_{a,n})R_2(x,\lambda_{a,n})
   +2
   R_1(x,\lambda_{a,n})\,\varepsilon_n^1+2R_2(x,\lambda_{a,n})\,\varepsilon_n^2\\ %
   &\phantom{=}+\varepsilon_n^1\varepsilon_n^2\Big)\n{\sr_n(\cdot,p,q)}{2}^{-2}.
\end{eqnarray*}
From \eqref{akns-eq_asymptotique_sol_reg_L2}, we have
$$\left|\varepsilon_n^j(x)\right|\leq \left(\frac{x}{1+|\lambda_{a,n}|x}\right)^a\ell^2(n),\quad j=1,2.$$
Thus, using \eqref{akns-est-norm-Rn}, we get
\begin{eqnarray*}
\fl\qquad  2G_{n,1}(x,V)G_{n,2}(x,V) =& -2j_a(\lambda_{a,n}x)j_{a-1}(\lambda_{a,n}x)(1+\ell^2(n))\\
   &+2{\lambda_{a,n}}^a \Big(j_{a-1}(\lambda_{a,n}x)\,\varepsilon_n^2-j_a(\lambda_{a,n}x)\,\varepsilon_n^1\Big)+\ell^1(n)
\end{eqnarray*}
and
\begin{eqnarray*}
 \fl\qquad G_{n,2}(x,V)^2-G_{n,1}(x,V)^2 =& \big(j_a(\lambda_{a,n}x)^2-j_{a-1}(\lambda_{a,n}x)^2\big)(1+\ell^2(n))\\
   &+2{\lambda_{a,n}}^a
   \Big(-j_{a-1}(\lambda_{a,n}x)\,\varepsilon_n^1-j_a(\lambda_{a,n}x)\,\varepsilon_n^2\Big)+\ell^1(n),
\end{eqnarray*}
then, we obtain uniformly for $x\in[0,1]$,
 \begin{eqnarray*}
  \fl\qquad r_n(x,V)=2{\lambda_{a,n}}^a&\Big[\ja{a-1}{_{a,n} x} \varepsilon_n(x,V)^\perp-\ja{a}{_{a,n}
x}
\varepsilon_n(x,V)\Big]\\&+\Phi_a(\lambda_{a,n}x)\ell^2(n)+\ell^1(n).
\end{eqnarray*}
With the uniform estimation on $[\delta,1]$, $\ja{a}{_{a,n}
x}=\sin\left(\lambda_{a,n}
x-\frac{a\pi}{2}\right)+\mathcal{O}\left(\frac{1}{\lambda_{a,n}}\right)$,
we get
\begin{eqnarray*}
  \left\langle\omega,e_{n,1}\right\rangle&=&\left\langle \omega,{B_a}^\ast[r_n]\right\rangle=\left\langle B_a[\omega],r_n\right\rangle\\
   &=& \int_0^1 \cos\left(\lambda_{a,n}
t-\frac{a\pi}{2}\right) 2{\lambda_{a,n}}^a\varepsilon_n(t,V)^\perp\cdot B_a[\omega](t)d t\\
    && -\int_0^1 \sin\left(\lambda_{a,n}
t-\frac{a\pi}{2}\right) 2{\lambda_{a,n}}^a\varepsilon_n(t,V)\cdot
B_a[\omega](t)d t\\
    &&+\left\langle \ell^2(n) B_a[\omega],\Phi_a(\lambda_{a,n}x)\right\rangle+\ell^1(n).
\end{eqnarray*}
Now, with lemma \ref{akns-lemme_AG_A1}, notice that for all
$f\in\lr$, we have uniformly on the bounded sets of $\lr$,
$$\left|\int_0^1 \cos\left(\lambda_{a,n} x\right) f(t) d t\right|=\n{f}{2}\left|\int_0^1
\cos\left(\lambda_{a,n} x\right) \frac{f(t)}{\n{f}{2}} d
t\right|\leq \n{f}{2}\,\ell^2(n).$$ This leads for instance to
\begin{eqnarray*}
\fl\qquad  \left|\int_0^1 \sin\left(\lambda_{a,n}
t-\frac{a\pi}{2}\right) 2{\lambda_{a,n}}^a\varepsilon_n(t,V)\cdot
B_a[\omega](t)d t\right| &\leq 2\ell^2(n)
\n{{\lambda_{a,n}}^a\varepsilon_n\cdot B_a[w]}{2},\\
&\leq 2 \ell^2(n)\ell^2(n)
\n{B_a[w]}{2},\\
&\leq\ell^1(n)\n{B_a[w]}{2}.
\end{eqnarray*}
And with the transformation operator, we get $\left\langle
\ell^2(n)
B_a[\omega],\Phi_a(\lambda_{a,n}x)\right\rangle=\ell^1(n)$.
Consequently, we have
$\left\langle\omega,e_{n,1}\right\rangle=\ell^1(n)$.

Now let $\Sigma_n=(\Sigma_n^1,\Sigma_n^2)$ be defined by
$\Sigma_n(x,V)=\ss_n(x,V)-S(x,\lambda_{a,n})$. With
\eqref{akns-eq_asymptotique_sol_sing_L2}, we have
$$\left|\Sigma_n^j(x)\right|\leq \left(\frac{1+|\lambda_{a,n}|x}{x}\right)^a\ell^2(n),\quad j=1,2.$$
First, with the definition of $A_n(x,p,q)$ and relations
\eqref{akns-eq_asymptotique_sol_reg_L2} and
\eqref{akns-eq_asymptotique_sol_sing_L2}, we have
\begin{eqnarray*}
\fl\qquad    A_n(x,p,q)=\Psi_a(\lambda_{a,n}x)
        &-{\lambda_{a,n}}^{-a}\big(j_{a-1}(\lambda_{a,n}x){\Sigma_n(x,V)}^\perp
        -j_a(\lambda_{a,n}x)\Sigma_n(x,V)\big)\\
        &+{\lambda_{a,n}}^a\big(\eta_{a-1}(\lambda_{a,n}x){\varepsilon_n(x,V)}^\perp
        -\eta_a(\lambda_{a,n}x)\varepsilon_n(x,V)\big)+\ell^1(n),
\end{eqnarray*}
which leads, using \eqref{akns-expr_grad_kappa} with
\eqref{akns-produit-scalaire}, to
\begin{eqnarray*}
\fl\qquad    \frac{\nabla_{p,q}\kappa_{a,n}}{\kappa_{a,n}}=&\Psi_a(\lambda_{a,n}x)+\ell^2(n)\Psi_a(\lambda_{a,n}x)\\
        &-{\lambda_{a,n}}^{-a}\left(j_{a-1}(\lambda_{a,n}x){\Sigma_n(x,V)}^\perp
        -j_a(\lambda_{a,n}x)\Sigma_n(x,V)\right)\\
        &+{\lambda_{a,n}}^a\left(\eta_{a-1}(\lambda_{a,n}x){\varepsilon_n(x,V)}^\perp
        -\eta_a(\lambda_{a,n}x)\varepsilon_n(x,V)\right)+\ell^1(n).
\end{eqnarray*}
Then, we get
\begin{eqnarray*}
\fl\qquad
s_n(x)=&-{\lambda_{a,n}}^{-a}\big(j_{a-1}(\lambda_{a,n}x){\Sigma_n(x,V)}^\perp
        -j_a(\lambda_{a,n}x)\Sigma_n(x,V)\big)\\
       & +{\lambda_{a,n}}^a\big(\eta_{a-1}(\lambda_{a,n}x){\varepsilon_n(x,V)}^\perp
        -\eta_a(\lambda_{a,n}x)\varepsilon_n(x,V)\big)\\
       & +\ell^2(n)\Psi_a(\lambda_{a,n}x)+\ell^2(n)\Psi_a(\lambda_{a,n}x)+\ell^1(n).
\end{eqnarray*}
Now, with the same arguments as previously, using the
transformation operator we find that
$$\left\{\left\langle\omega,e_{n,2}\right\rangle\right\}_{n\in\Z}\in\ell^1(\Z).$$
Thus, proposition \ref{prop-base} proves the result. \qed

We can go further in solving the inverse spectral problem. Indeed,
we can give explicitly the inverse of the spectral map's
differential. But first some notations:

\begin{nota} For all $n\in\Z$, we set%
$$X_{a,n}(p,q)=\frac{-{\nabla_{p,q} \kappa_{a,n}}^\perp}{\kappa_{a,n}(p,q)},\quad
Y_{a,n}(p,q)=\frac{(-1)^n\,{\nabla_{p,q}
\lambda_{a,n}}^\perp}{\left[\left(|n|+\frac{a}{2}\right)\pi\right]^a\kappa_{a,n}(p,q)}.$$
\end{nota}

Notice that, according to estimations from corollary
\ref{akns-thm_asympt_misc}, we have
\begin{equation}\label{akns-est-Xn-Yn}
    X_{a,n}(p,q)=-{\Psi_a{(\lambda_{a,n}x)}}^\perp+\ell^2(n),\, Y_{a,n}(p,q)={\Phi_a{(\lambda_{a,n}x)}}^\perp+\ell^2(n).
\end{equation}

\begin{coro}\label{akns-exp-inverse-diff_akns} $\lambda^a\times\kappa^a$ is a local
real analytic diffeomorphism at every point in $\lr\times\lr$% dans $\ell^2(\Z)\times\ell^2(\Z)$
. Moreover, the inverse of
$d_{p,q}\left(\lambda^a\times\kappa^a\right)$ is the linear map
from $\ell^2_\R(\Z)\times\ell^2_\R(\Z)$ onto
$\lr\times\lr$ given by%
$$\left(d_{p,q}\left(\lambda^a\times\kappa^a\right)\right)^{-1}(\xi,\eta)=\sum_{n\in\Z}\xi_n X_{a,n}+\sum_{n\in\Z}\eta_n Y_{a,n}.$$
\end{coro}

\noindent\textit{Proof.} First point comes directly from the
theorem and the definition of a local diffeomorphism. Now consider
$(\xi,\eta)\in\ell^2_\R(\Z)\times\ell^2_\R(\Z)$ and let%
$$u=\sum_{n\in\Z}\xi_n X_{a,n}+\sum_{n\in\Z}\eta_n Y_{a,n}.$$
Thanks to relation \eqref{akns-akns_eq_lemme_bessel_sinus_3}, the
transformation operator lets us write estimations
\eqref{akns-est-Xn-Yn} in the following way
$$X_{a,n}(p,q)={B_a} \Bigg[\Bigg[\begin{array}{c}
  \sin(2\lambda_{a,n}x) \\
  \cos(2\lambda_{a,n}x) \\
\end{array}\Bigg]+\ell^2(n)\Bigg],%$$
%$$
Y_{a,n}(p,q)={B_a} \Bigg[\Bigg[\begin{array}{c}
  \cos(2\lambda_{a,n}x) \\
  -\sin(2\lambda_{a,n}x) \\
\end{array}\Bigg]+\ell^2(n)\Bigg].$$ Since $B_a$ is bounded and
$\xi$, $\eta$ are in $\ell^2_{\R}(\Z)$, the sum defining $u$
exists in $\lr\times\lr$ . Orthogonality relations from section
\ref{akns-section-ortho} imply that for
all $n\in\Z$%
$$\left\langle \nabla_{p,q}\lambda_{a,n},u\right\rangle=\xi_n \quad
\mathrm{ et } \quad \left\langle \nabla_{p,q}\widetilde{\kappa}_{a,n},u\right\rangle=\eta_n.$$%
Thus we have
$d_{p,q}\left(\lambda^a\times\kappa^a\right)(u)=(\xi,\eta)$, which
proves the corollary. \qed

%\subsection{Isospectral sets}

We finish the local inverse spectral problem with the description
of isospectral sets. For $(p_0,q_0)\in\lr\times\lr$, we define the
set of AKNS potentials with same spectrum as $(p_0,q_0)$, called
isospectral set of $(p_0,q_0)$, by:
$$\mathrm{Iso}(p_0,q_0,a)=\left\{(p,q)\in\lr\times\lr:
\lambda^a(p,q)=\lambda^a(p_0,q_0)\right\}.$$

\begin{thm} Let $(p_0,q_0)\in\lr\times\lr$, then%
    \begin{enumerate}[(a)]
        \item $\mathrm{Iso}(p_0,q_0,a)$ is a real analytic submanifold of $\lr\times\lr$.
        \item At every point $(p,q)$ of
        $\mathrm{Iso}(p_0,q_0,a)$, the tangent space is
        $$T_{p,q}\mathrm{Iso}(p_0,q_0,a)=\left\{\sum_{n\in\Z}\eta_n Y_{a,n}(p,q)\,: \eta\in\ell^2_\R(\Z)\right\}$$
    and the normal space is
        $$N_{p,q}\mathrm{Iso}(p_0,q_0,a)=\left\{\sum_{n\in\Z}\eta_n {Y_{a,n}(p,q)}^\perp\,: \eta\in\ell^2_\R(\Z)\right\}.$$
\end{enumerate}
\end{thm}

\noindent\textit{Proof.}    %Let $(p,q)\in\mathrm{Iso}(p_0,q_0,a)$.
    %According to
    %corollary \ref{akns-exp-inverse-diff_akns}, $\lambda^a\times\kappa^a$
    %is a real-analytic diffeomorphism locally in the neighborhood
    %of $(p,q)$. Thus there exist $\mathcal{V}$ open neighborhood of $(p,q)$ in $\lr\times\lr$ and
    %an open subset $\mathcal{W}$ of $\ell^2_\R(\Z)\times\ell^2_\R(\Z)$ such that
    %$\lambda^a\times\kappa^a: \mathcal{V}\rightarrow\mathcal{W}$
    %is a real analytic isomorphism. We then have
    %$$\left(\lambda^a\times\kappa^a\right)(\mathcal{V}\cap\mathrm{Iso}(p_0,q_0,a))=
    %\mathcal{W}\cap\left(\left\{\lambda^a(p_0,q_0)\right\}\times \ell^2_\R(\Z)\right),$$
    %which proves point \textit{(a)}.\\
    Notice that the local real-analytic diffeomorphism $\lambda^a\times\kappa^a$ defines a
    chart at each point $(p,q)\in\mathrm{Iso}(p_0,q_0,a)$, the definition of a
    submanifold gives point \textit{(a)}.\\
%    \item
    Since $T_{p,q}\mathrm{Iso}(p_0,q_0,a)=
    \left(d_{p,q}\left(\lambda^a\times\kappa^a\right)\right)^{-1}\big(\big\{0_{\ell^2_\R(\Z)}\big\}\times\ell^2_\R(\Z)\big)$,
    corollary \ref{akns-exp-inverse-diff_akns} gives the expression of the tangent space.
    Now, the family $(Y_{a,n})_{n\in\Z}$ is free since $(\nabla_{p,q}\lambda_{a,n})_{n\in\Z}$ is. Moreover, it is orthogonal to
    $({Y_{a,n}}^\perp)_{n\in\Z}$. Then we have the first inclusion
    $$\left\{\sum_{n\in\Z}\eta_n {Y_{a,n}(p,q)}^\perp\,: \eta\in\ell^2_\R(\Z)\right\}\subset N_{p,q}\mathrm{Iso}(p_0,q_0,a).$$
    Now, every vector orthogonal to $({Y_{a,n}}^\perp)_{n\in\Z}$ is orthogonal to the gradients $(\nabla_{p,q}\lambda_{a,n})_{n\in\Z}$,
    in other words, is in the kernel of $d_{p,q}\lambda^a$. Thus the second inclusion follows and so does point \textit{(b)}. \qed

\subsection{A Borg-Levinson theorem on $\hunr\times\hunr$}

\begin{thm}\label{akns-th_injectif_borg} $\lambda^a\times\kappa^a$ is one-to-one
on $\hunr\times\hunr$.
\end{thm}

As in the case of a radial Schrödinger operator (see for instance
\cite{carlson}),  we introduce another solution to
\eqref{akns-AKNS+V} with boundary condition at $x=1$.

\begin{lem}\label{akns-lemme_estimation_rho}
Let $\rho(x,\lambda,V)$ be the solution of \eqref{akns-AKNS+V}
such that
\begin{equation}\label{akns-condition_en_1}
    \rho(1,\lambda,V)={u_\beta}^\perp.%\Bigg[\begin{array}{c}
%      \cos{\beta} \\
%      -\sin{\beta} \\
%    end{array}\Bigg].
\end{equation}
Then $\rho$ verifies the following properties
\begin{enumerate}[(i)]
    \item For $V=(p,q)\in\lc\times\lc$ and $\delta>0$, uniformly on $[\delta,1]$,
    $$\left|\rho(x,\lambda,V)-\Bigg[\begin{array}{c}
      \cos{\left(\lambda(1-x)-\beta\right)} \\
      \sin{\left(\lambda(1-x)-\beta\right)} \\
    \end{array}\Bigg]\right|\leq K(x) e^{|\im\lambda|(1-x)}$$
    where
    $K(x)=\displaystyle\exp{\left[\int_x^1 \left(|p(t)|+|q(t)|+\frac{a}{t}\right)d
    t\right]}$.
    \item For $V=(p,q)\in H^1\times H^1$ and $\delta>0$, uniformly on $[\delta,1]$,
$$\left|\rho(x,\lambda,V)-\Bigg[\begin{array}{c}
      \cos{\left(\lambda(1-x)-\beta\right)} \\
      \sin{\left(\lambda(1-x)-\beta\right)} \\
    \end{array}\Bigg]\right|\leq C_a\frac{K(x)}{\lambda x} (\n{V}{H^1}+1)e^{|\im\lambda|(1-x)}$$
    \item For all $x\in(0,1]$, $\rho(x,\lambda,V)$ is analytic on $\C\times\lc\times\lc$.
    \item For $n\in\Z$ and $\lambda=\lambda_{a,n}(V)$, we have
    \begin{equation}\label{akns-lien_R_rho}
    \sr_n(x,V)=\kappa_{a,n}(V)\rho(x,\lambda_{a,n}(V),V).
    \end{equation}
\end{enumerate}
\end{lem}

\noindent \textit{Lemma's proof.}$\,$\\ Points \textit{(i)},
\textit{(ii)} et \textit{(iii)} follow directly from a Picard
iteration construction of $\rho$. Indeed, we define as in the
regular case (see for instance \cite{gg}) $\rho$ with

Now prove point \textit{(iv)}: when $\lambda=\lambda_{a,n}(V)$,
according to \eqref{akns-AKNS-bord} and
\eqref{akns-condition_en_1}, $\rho(1,\lambda_{a,n}(V),V)$ and
$\sr(1,\lambda_{a,n}(V))$ are collinear. Then
$\rho(x,\lambda_{a,n}(V),V)$ and $\sr(x,\lambda_{a,n}(V))$
solutions of \eqref{akns-AKNS+V} with the same eigenvalue
$\lambda$ are also collinear, in other words there exists
$C_n\in\R$ such that $\sr_n(x,V)=C_n\,\rho(x,\lambda_{a,n}(V),V)$.
Using again \eqref{akns-AKNS-bord} then
\eqref{akns-condition_en_1} and \eqref{akns-def_kappan}, we deduce
that $\kappa_{a,n}(V)=C_n$. \qed

\noindent\textit{Proof of Theorem \ref{akns-th_injectif_borg}.}
Let $V,W\in\lr\times\lr$ such that
$(\lambda^a\times\kappa^a)(V)=(\lambda^a\times\kappa^a)(W)$. %, in
%other words
%
%$$ \forall n\in\Z,\quad
%  \lambda_{a,n}(V)= \lambda_{a,n}(W),\quad
%  \kappa_{a,n}(V)= \kappa_{a,n}(W).$$
For $u\in\R^2$, introduce the function
$$f(x,\lambda,V,W)=\frac{\left[\sr(x,\lambda,V)\negthinspace\cdot\negthinspace u-
\sr(x,\lambda,W)\negthinspace\cdot\negthinspace
u\right]\left[\rho(x,\lambda,V)\cdot u-\rho(x,\lambda,W)\cdot
u\right]}{D(\lambda,V)}.$$ For all $x\in(0,1]$, $f:\lambda\mapsto
f(x,\lambda,V,W)$ is a meromorphic function on $\C$ which has
simple poles $\lambda_{a,n}(V)$, $n\in\Z$. From the simplicity of
poles and since $f(\lambda)=h(\lambda)/g(\lambda)$, the residue of
$f$ at $\lambda_{a,n}(V)$ is
$$\mathrm{Res}(f,\lambda_{a,n}(V))=\frac{h(\lambda_{a,n}(V))}{g'(\lambda_{a,n}(V))}.$$

%First, with
Using that $\lambda_{a,n}(V)=\lambda_{a,n}(W)$ %we have
%$$\sr(x,\lambda_{a,n}(V),W)=\sr(x,\lambda_{a,n}(W),W)=\sr_n(x,W),$$
%then with relation \eqref{akns-lien_R_rho} from the preceding
%lemma, we get
%$$\rho(x,\lambda_{a,n}(V),V)=\frac{1}{\kappa_{a,n}(V)}\sr_n(x,V),$$
%and
%$$\rho(x,\lambda_{a,n}(V),W)=\rho(x,\lambda_{a,n}(W),W)=\frac{1}{\kappa_{a,n}(W)}\sr_n(x,W).$$
%Since
and $\kappa_{a,n}(V)=\kappa_{a,n}(W)$, together with relations
\eqref{akns-lien_R_rho} and
%
%
%we deduce that
%$$\mathrm{Res}(f,\lambda_{a,n}(V))=\frac{\left[\sr_n(x,V)\cdot u-
%\sr_n(x,W)\cdot u\right]^2}{\kappa_{a,n}(V) \frac{\partial
%D}{\partial \lambda}(\lambda_{a,n}(V),V)}$$ and, with
\eqref{akns-eq_spectre_simple}, we obtain
$$\mathrm{Res}(f,\lambda_{a,n}(V))=-\frac{\left[\sr_n(x,V)\cdot u-\sr_n(x,W)\cdot
u\right]^2}{\n{\sr_n(\cdot,V)}{2}^2}.$$ %
To conclude, we make use of a complex analysis result
\begin{lem}[Lemma 3.2 \cite{ist}]\label{akns-lemme_residus_nuls} Let $f$ be a meromorphic function on $\C$ such
that
$$\sup_{|\lambda|=r_n}|f(\lambda)|=o\left(\frac{1}{r_n}\right)$$
for an unbounded sequence of positive real numbers $(r_n)$. Then,
the sum of the residues of $f$ is zero.
\end{lem}

Let $N>0$ be an integer and $C_N$ be the circle defined by
$$\left|\lambda-\left(\frac{a\pi}{2}+\beta\right)\right|=\left(N+\frac{1}{2}\right)\pi.$$
Estimate $|\lambda f(x,\lambda,V,W)|$ on $C_N$. From
\eqref{akns-estimation_h1_reg} and \eqref{akns-estimation_L2_loc}
with the help of lemma \eqref{akns-lemme_estimation_rho}, we have
for $N$ large enough
\begin{eqnarray*}
\fl\qquad  \left|\sr(x,\lambda,V)\cdot u-
\sr(x,\lambda,W)\cdot u\right| \leq C(\n{V}{H^1}+\n{W}{H^1})  e^{|\im\lambda|x}\frac{\ln{|\lambda|}}{|\lambda|^{a+1}}, \\
\fl\qquad   \left|\rho(x,\lambda,V)\cdot u-\rho(x,\lambda,W)\cdot
u\right|\leq\frac{K(x)}{|\lambda| x}
(\n{V}{H^1}+\n{W}{H^1})e^{|\im\lambda|(1-x)},\\
\fl\qquad|D(\lambda,V)|\geq|R(1,\lambda)\cdot
u_\beta|-|(\sr(1,\lambda,V)-R(1,\lambda))\cdot u_\beta|\geq
\frac{C}{|\lambda|^a}e^{|\im\lambda|}.
\end{eqnarray*}
We deduce that uniformly for $x\in[\delta,1]$ and $\lambda\in
C_N$, $\displaystyle|\lambda f(\lambda,V,W)|\leq C
\frac{\ln{|\lambda|}}{|\lambda|}$. Thus, %Ainsi, pour tout
result from lemma \ref{akns-lemme_residus_nuls} is valid  for $f$.
Since residues of $f$ have the same sign, they are all zero. In
conclusion, we have for all $n\in\Z$, $u\in\R^2$, $\delta\in(0,1]$
and $x\in[\delta,1]$, $\sr_n(x,V)\cdot u-\sr_n(x,W)\cdot u=0$. We
can deduce, recalling continuousness of eigenvectors at $x=0$,
that for all $x\in [0,1]$ and all $n\in\Z$
$$\sr_n(x,V)=\sr_n(x,W).$$
Plug this in \eqref{akns-AKNS+V} to deduce that % , nous en
%déduisons que
%$$(V(x)-W(x))\sr_n(x,V)=0,\quad n\in\Z.$$
$V=W$ almost every where on $[0,1]$. \qed

\ack The author would like to thank B. Grébert and J.-C. Guillot
for numerous interesting discutions, and to acknowledge the
hospitality of the Institüt fur Mathematik of the University of
Zürich, namely T. Kappeler, where the present work was completed.

%   Annexes
    %   Annexe A: Fonctions de Bessel
\appendix \section*{Appendix}\setcounter{section}{1}\label{annexe-bessel}

Spherical Bessel functions $j_a$ and $\eta_a$ are defined through
\begin{equation}\label{annexe-bessel-def}
    j_a(z)=\sqrt{\frac{\pi z}{2}}J_{a+1/2}(z),\quad
    \eta_a(z)=(-1)^a \sqrt{\frac{\pi z}{2}}J_{-a-1/2}(z),
\end{equation}
where $J_\nu$ is the first kind Bessel function of order $\nu$
(see \cite{emot} for precisions).

The following estimates can be found in \cite{moi-these}.
\begin{itemize}
    \item Uniform estimates on $\C$:
    \begin{eqnarray}
      \label{annexe-bessel_est_ja}\left|j_{a}(z)\right| &\leq& C e^{|\im z|}\left(\frac{|z|}{1+|z|}\right)^{a+1}, \\
      \label{annexe-bessel_est_na}\left|\eta_{a}(z)\right| &\leq& C e^{|\im z|}\left(\frac{1+|z|}{|z|}\right)^{a}.
    \end{eqnarray}
    \item Estimations for the Green function $G(x,t,\lambda)$ when $0\leq t\leq x$:
    \begin{equation}\label{annexe-resolvante_est_tx}%
        |\G(x,t,\lambda)| \leq C e^{|\im
        \lambda|(x-t)}\left(\frac{x}{1+|\lambda|x}\right)^{a}\left(\frac{1+|\lambda|t}{t}\right)^{a}.
    \end{equation}
    \item Estimations for the Green function $G(x,t,\lambda)$ when $0\leq x\leq t\leq
    1$:
    \begin{equation}\label{annexe-resolvante_est_xt}%
        |\G(x,t,\lambda)|\leq C e^{|\im
        \lambda|(t-x)}\left(\frac{1+|\lambda|x}{x}\right)^{a}\left(\frac{t}{1+|\lambda|t}\right)^{a}.
    \end{equation}
    \item   Trigonometric expression (\cite{emot} formulas $(1-2)$ section 7.11 p.78),
    \begin{eqnarray}\label{annexe-bessel-ja-sinus}
    j_a(z)=\sin{\left(z-\frac{a\pi}{2}\right)}P_a(z^{-1})+\cos{\left(z-\frac{a
        \pi}{2}\right)}I_a(z^{-1}),\\
%    \end{equation}
%    \begin{equation}
    \label{annexe-bessel-na-sinus}
    \eta_{a}(z)=\cos{\left(z-\frac{a\pi}{2}\right)}P_a(z^{-1})-\sin{\left(z-\frac{a\pi}{2}\right)}I_a(z^{-1})
    \end{eqnarray}
    where $P_a$ and $I_a$ are even, resp. odd, polynomials given
    by
    \begin{eqnarray}\label{annexe-bessel-Pa}
    \fl\qquad    P_a(z)=\sum_{m=0}^{\leq
    a/2}(-1)^m
    \big(a+1/2,2m\big)(2z)^{2m},\phantom{+1^{+1}}\quad&(P_a(0)=1),\\
    \label{annexe-bessel-Ia}
    \fl\qquad    I_a(z)=\sum_{m=0}^{\leq
    (a-1)/2}(-1)^m
    \big(a+1/2,2m+1\big)(2z)^{2m+1},\quad&(I_a(0)=0),
    \end{eqnarray}
    where %$(\nu,m)$
$\displaystyle
%    \begin{equation}\label{annexe-bessel-hankel}
    (\nu,m)=\frac{\Gamma(\nu+1/2+m)}{m!\Gamma(\nu+1/2-m)}
%\end{equation}
%
$ is the Hankel symbol.
\end{itemize}

%%%%%%%%%%%%%%%%%%%%%%%%%%%%%%%%%%%%%%%%%%%%%%%%%%%%%%%%%%%%%%%%%%%%%%%%%%%%%%%%%%%%%%%%%%%%%%
\subsection{Technical lemmas}

%Tout d'abord, nous aurons besoin de deux lemmes techniques:
%%%%%%%%%%%%%%%%%%%%%%%%%%%%%%%%%%%%%%%%%%%%%%%%%%%%%%%%%%%%%%%%%%%%%%%%%%%%%%%%%%%%%%%%%%%%%%
\begin{lema}\label{annexe-lemme_primitive_K} Let $f_1(z)=2 j_{a-1}(z)j_{a}(z)$. Then $F_1=\int f_1(z)dz$ such that
$F_1(0)=0$ verifies  the properties
\begin{enumerate}[(i)]
    \item $\displaystyle\left|F_1(z)\right|\leq C
    \left(\frac{|z|}{1+|z|}\right)^{2a+2}$ for $|z|\leq 1$;
    \item $\displaystyle F_1(z)=-a \mathrm{ci}(2z)+p_a\left(z^{-1}\right)\cos{(2z)}+
    q_a\left(z^{-1}\right)\sin{(2z)}+r_a\left(z^{-1}\right)$
    if $z\neq 0$,
\end{enumerate} Where $\displaystyle\mathrm{ci}(z)=\int_0^z
\frac{\cos{t}-1}{t}\,d t$ and $p_a$, $q_a$, $r_a$ are resp. even,
odd and even, polynomials.
\end{lema}
%%%%%%%%%%%%%%%%%%%%%%%%%%%%%%%%%%%%%%%%%%%%%%%%%%%%%%%%%%%%%%%%%%%%%%%%%%%%%%%%%%%%%%%%%%%%%
\begin{lema}\label{annexe-lemme_primitive_L} Let $f_2(z)=\eta_{a-1}(z)j_{a}(z)+\eta_{a}(z)j_{a-1}(z)$. Then
$F_2=\int f_2(z)dz$ such that $F_2(0)=0$ satisfies the properties
\begin{enumerate}[(i)]
    \item $\displaystyle\left|F_2(z)\right|\leq C
    \frac{|z|}{1+|z|}$ for $|z|\leq 1$;
    \item $\displaystyle F_2(z)=a \mathrm{Si}(2z)-p_a\left(z^{-1}\right)\sin{(2z)}+
    q_a\left(z^{-1}\right)\cos{(2z)}$ if $z\neq 0$.
\end{enumerate}
Where $\displaystyle\mathrm{Si}(z)=\int_0^z \frac{\sin{t}}{t}\,d
t$ and $p_a$, $q_a$ are the previous polynomials.
\end{lema}

    \subsection{Calculation lemma} The following lemma is adapted from
\cite{carlson}, its proof lies on some Hardy inequalities (for
details see \cite{carlson} and \cite{moi-these}). Together with
the transformation operator, it is an essential tool for the
computation of asymptotics for $\lr$ potentials.

%%%%%%%%%%%%%%%%%%%%%%%%%%%%%%%%%%%%%%%%%%%%%%%%%%%%%%%%%%%%%%%%%%%%%%%%%%
\begin{lema}[Carlson \cite{carlson}]\label{annexe-lemme-carlson-1} Let $f\in\lc$ and $(z_n)_{n\in\N}$ a strictly positive real sequence
        such that
        \begin{equation*}%\label{annexe-eq-lemme-carlson-2}
           z_0>0\quad\mathrm{ and }\quad \exists(C_1,C_2)\in\R^\ast_+\times\R^\ast_+,\,\forall n\in\N,\quad C_1\leq z_{n+1}-z_n\leq C_2.
        \end{equation*}
        Then, uniformly on bounded set in $\lc$,
        $$\left(\int_0^{1/z_n}|f(t)|dt\right)_{n\in\N}, \left(\int_{1/z_n}^1\left|\frac{f(t)}{z_n t}\right|dt\right)_{n\in\N}\in\ell^2_\R(\N).$$
\end{lema}
%%%%%%%%%%%%%%%%%%%%%%%%%%%%%%%%%%%%%%%%%%%%%%%%%%%%%%%%%%%%%%%%%%%%%%%%%%

%   Bibliographie
%GATHER{Xbib.bib}   % For Gather Purpose Only
\section*{References}
\bibliographystyle{plain}
\bibliography{Xbib}
\end{document}